\documentclass[mathpazo]{csiam-am}
\renewcommand\thisnumber{x}
\renewcommand\thisyear {20xx}
\renewcommand\thismonth{xx}
\renewcommand\thisvolume{x}

\usepackage{xcolor,txfonts}
\usepackage{bm}
\usepackage{algorithm,algpseudocode, subfigure, multirow}
\usepackage{booktabs,makecell}
\usepackage{graphicx}

\begin{document}
\title[pETNNs]{pETNNs: Partial Evolutionary Tensor Neural Networks for Solving Time-dependent Partial Differential Equations}


 \author[Kao T. et.~al.]{Tunan Kao\affil{1}, He Zhang\affil{1}, Lei Zhang\affil{1, 2, \corrauth}, Jin Zhao\affil{3, 4}}
 \address{\affilnum{1}\ {Beijing International Center for Mathematical Research, Peking University, Beijing, 100871,  China}\\
\affilnum{2}\ {Center for Quantitative Biology, Center for Machine Learning Research, Peking University, Beijing, 100871, China}\\
\affilnum{3}\ {Academy for Multidisciplinary Studies, Capital Normal University, Beijing, 100048, China}\\
\affilnum{4}\ {Beijing National Center for Applied Mathematics, Beijing, 100048, China}\\
 }



%
%
 \emails{{\tt kaotunan@pku.edu.cn} (T.~Kao), {\tt zhanghe@bicmr.pku.edu.cn} (H.~Zhang), {\tt pkuzhangl@pku.edu.cn} (L.~Zhang), {\tt zjin@cnu.edu.cn} (J.~Zhao) }

\begin{abstract}
We present the partial evolutionary tensor neural networks (pETNNs), a novel framework for solving time-dependent partial differential equations with high accuracy and capable of handling high-dimensional problems. Our architecture incorporates tensor neural networks and evolutionary parametric approximation. A posterior error bounded is proposed to support the extrapolation capabilities. In the numerical implementations, we adopt a partial update strategy to achieve a significant reduction in computational cost while maintaining precision and robustness. Notably, as a low-rank approximation method of complex dynamical systems, pETNNs enhance the accuracy of evolutionary deep neural networks and empower computational abilities to address high-dimensional problems. Numerical experiments demonstrate the superior performance of the pETNNs in solving time-dependent complex equations, including the incompressible Navier-Stokes equations, high-dimensional heat equations, high-dimensional transport equations, and dispersive equations of higher-order derivatives.
\end{abstract}

\ams{65M99, 65L05, 68T07
}
\keywords{time-dependent partial differential equations, tensor neural networks, evolutionary deep neural networks, high-dimensional problems.}

\maketitle

\section{Introduction}
Partial differential equations (PDEs) are ubiquitous in modeling phenomena across scientific and engineering disciplines. They serve as indispensable tools in modeling continuum mechanics, electromagnetic theory, quantum mechanics, and a myriad of other fields where the evolution of systems across space and time is of interest. Traditional numerical approaches for solving PDEs, such as finite difference \cite{leveque2007finite}, finite element \cite{ciarlet2002finite}, and spectral methods \cite{shen2011spectral}, have been widely used. However, the computational burden imposed by these methods grows exponentially with the increase in the problem dimensionality, often rendering them impractical for high-dimensional systems. This phenomenon, known as the ``curse of dimensionality", has been a persistent impediment to progress in various scientific domains.

The emergence of machine learning has introduced a novel set of tools to the scientific community, offering a potential panacea to the curse of dimensionality. Deep learning, a class of machine learning characterized by deep neural networks (DNNs), has been particularly successful in areas where traditional algorithms falter due to the complexity and volume of the data involved, such as \cite{cnm.1640100303, Lagaris712178, long2018pde1, Geist2021, Jin2023, Zhang2024}. The universal approximation theorem underpins this capability, suggesting that a neural network can approximate any continuous function to a desired degree of accuracy \cite{Cybenko1989, Hornik1989}. Leveraging this, researchers have proposed various frameworks where DNNs are trained to satisfy the differential operators, initial conditions, and boundary conditions of PDEs.

A notable advancement in the field is the emergence of deep Galerkin method \cite{SIRIGNANO20181339}, deep Ritz method \cite{bing2018DRM}, and physics-informed neural networks (PINNs) \cite{RAISSI2019686}. They embed the governing physical laws, encapsulated by PDEs, into the architecture of deep learning models. By incorporating the PDEs directly into the loss function, PINNs ensure that the learned solutions are not merely data-driven but also conform to the underlying physical principles. This integration of physical laws into the learning process imbues PINNs with the ability to generalize beyond the data they were trained on, making them particularly adept at handling scenarios where data is scarce or expensive to acquire.

However, the efficacy of PINNs is predominantly confined to the temporal domain for which they have been trained, typically within the interval $[0,T]$. Their ability to extrapolate beyond this training window is limited, which is a manifestation of neural networks' inherent weakness in out-of-distribution generalization. This limitation hinders their predictive capacity, rendering them less effective for forecasting future states of the system under study.

The evolutionary deep neural networks (EDNNs) \cite{PhysRevE.104.045303}, which can address this challenge, have been developed as an innovative approach to solving time-dependent PDEs. The EDNNs are designed to evolve in tandem with the temporal dynamics they model, thus possessing an enhanced capability for prediction. This is achieved by structuring the neural network in a way that intrinsically accounts for the temporal evolution, allowing for a more robust extrapolation into future times. The methodology derived by the EDNNs has attracted significant attention. {  Hao et al. in \cite{hao2023NED} proposed a neural energy descent method, which identifies steady-state solutions of evolutionary equations to optimize neural networks. The work in \cite{NCPS2024} formulated deep neural network parameters as an optimal control problem to approximate solution operators of evolutionary PDEs.} The authors in \cite{BRUNA2024112588} employed this to build upon the foundational results established in \cite{doi:10.1137/21M1415972}. Specifically, they move away from the traditional method of training DNNs for PDEs and adopt a new strategy that uses the Dirac-Frenkel variational principle to evolve the network parameters through a system of ordinary differential equations (ODEs). Subsequently, they have enhanced the computational efficiency by randomized sparse neural Galerkin schemes in \cite{NEURIPS2023_0cb310ed}. Furthermore, the authors in \cite{kast2023positional} conducted a comprehensive investigation into the boundary treatment of EDNNs, leading to significant improvements.

Despite the method that the EDNNs can effectively compute the time-dependent PDEs and reasonably predict the related solutions, it seems unable to yield high precision due to its reliance on the Monte-Carlo approach in the process of computing integrals. On the other hand, many works that delve into various techniques aimed at improving the accuracy and efficiency of deep learning-based solvers for PDEs, such as \cite{JAGTAP2020109136,doi:10.1137/19M1274067, NEURIPS2022_374050dc, doi:10.1137/22M1527763, Muller2023achieving}. Among them, the tensor neural networks (TNNs) emerge as a methodological innovation \cite{liao2022solving,wang2023tensor,wang2023tensorSolving}, characterized by a restructured network architecture that can employ the Gaussian quadrature formula as an alternative to Monte Carlo integration, thus substantially enhance the accuracy of solutions. Nonetheless, the current implementation of this approach is exclusively tailored to stationary partial differential equations, and its utilization in the context of time-dependent differential equations has not been studied in the existing scholarly discourse.

This paper presents the partial evolutionary tensor neural networks (pETNNs, pronounced ``$Peten$''), enhancing the TNNs framework to robustly compute time-dependent equations with exceptional accuracy. Our contributions can be delineated as follows:
\begin{itemize}
 \item \textbf{Enhanced Computational Framework}: We have enhanced the TNNs framework to effectively handle time-dependent equations while maintaining high accuracy.

\item \textbf{Innovative Update Strategies}: Novel parameter update strategies have been introduced, where only partial parameters are updated. This strategic limitation optimizes the allocation of computational resources, maintaining the system's efficiency.

\item \textbf{Posterior Error Estimation}: We formulate the posterior error bound for pETNNs based on the dynamical parametric approximation theory, which holds for various spatial discretization schemes.

\item \textbf{Empirical Validation}: Numerical experiments with pETNNs reveal significant reductions in computational costs without compromising accuracy. The enhanced models not only achieve higher fidelity in numerical solutions than conventional EDNNs but also demonstrate superior capability in resolving complex equations such as the Navier-Stokes equations, simplified Korteweg-de Vries equations, high-dimensional heat equations, and high-dimensional transport equations.
\end{itemize}

This paper is organized as follows. We propose the main results of the pETNNs in Section 2, including the derivation process, and the corresponding algorithms. Extensive numerical experiments are carried out to validate the efficiency of the proposed pETNNs in Section 3. We finally present our conclusion and discussion in Section 4.

\section{Partial Evolutionary Tensor Neural Networks}\label{sec:2}

In this section, we will show the structure and derivation of pETNNs. We first review the architecture of TNNs, which have been used for solving stationary PDEs \cite{wang2023tensor}.

\subsection{Tensor Neural Networks}
\label{sec:TNN}

\begin{figure}[ht!]
  \centering
  \includegraphics[width=0.80\textwidth]{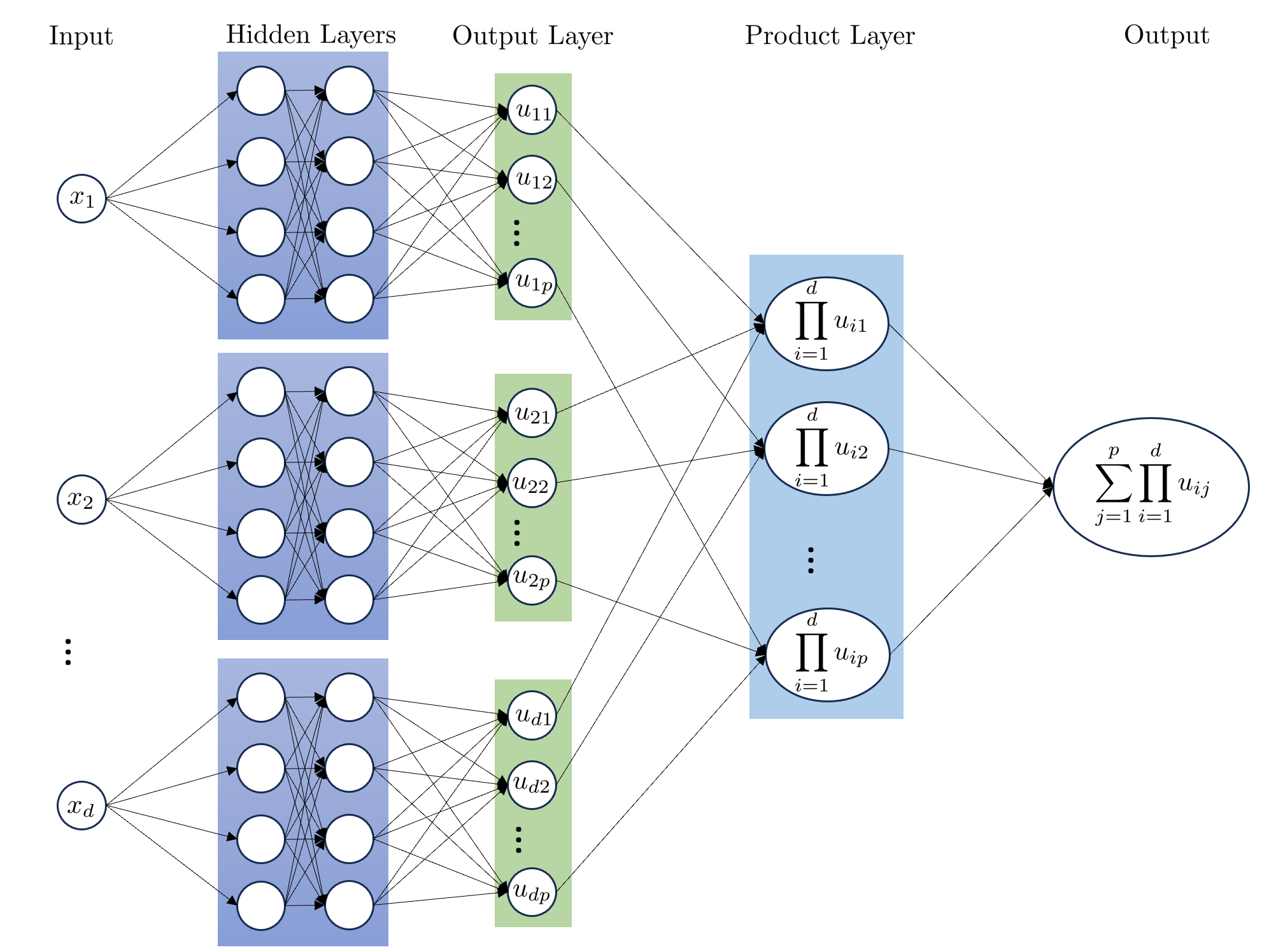}
  \caption{Schematic diagram of the tensor neural networks (TNNs). If the hidden layers of each sub-network is of width $M$ and depth $N$, then the scale of the total free parameters is $O(M^2Nd + Mpd)$, which grows linearly in terms of the problem dimension $d$.}
  \label{fig1}
\end{figure}

For simplicity, we primarily discuss the method of approximating a scalar function using the TNNs, and it should be noted that the approach can be extended to high-dimensional cases straightforwardly. The function $u(x)$ is approximated by the TNNs as follows:
\begin{equation}\label{3.1.1}
  u(x ; \theta)=\sum_{j=1}^p u_{1, j}\left(x_1 ; \theta^1\right) u_{2, j}\left(x_2 ; \theta^2\right) \cdots u_{d, j}\left(x_d ; \theta^d\right)=\sum_{j=1}^p \prod_{i=1}^d u_{i, j}\left(x_i ; \theta^i\right).
\end{equation}
Fig. \ref{fig1} shows the details about the TNNs with $d$ sub-networks. Here $\theta = \{\theta^1, \cdots, \theta^d\}$ denotes all parameters of the whole architecture, $d$ is the dimension, the positive integer $p$ is a hyper-parameter, which is related to the canonical polyadic (CP) decomposition\cite{2009Tensor} structure approximating to $u(x)$, and $u_{i, j}(x_i; \theta^i)$ is the $j$-th output of the $i$-th sub-network.

The introduction of TNNs drew inspiration from the MIONet \cite{MIONet2022}, leveraging its unique structure to enable highly precise and efficient numerical integration for solving high-dimensional PDEs. 
Specifically, one of the most advantages of TNNs is that the high-dimensional integration of functions represented by TNNs can be decomposed into the product of a series of one-dimensional integration, which can be calculated using the conventional Gaussian quadrature schemes and offers higher precision compared to the Monte Carlo methods. 
This approach can be effectively incorporated into the formulation of loss functions, especially those based on the Mean Squared Error (MSE) criterion, thereby improving the precision of solutions to PDEs obtained via deep neural network frameworks.
Furthermore, representing a function with the TNNs leads to a remarkable decrease in computational complexity for its numerical integration, which overcomes the curse of dimensionality in some sense (see Theorem 3 in \cite{wang2023tensor}). For the comprehensive discussion and the capabilities of the TNNs, we refer the readers to \cite{wang2023tensor,wang2023tensorSolving}.

\subsection{Partial Evolutionary Tensor Neural Networks}\label{sec:pETNN}

\begin{figure}[ht!]
  \centering
  \includegraphics[width=0.80\textwidth]{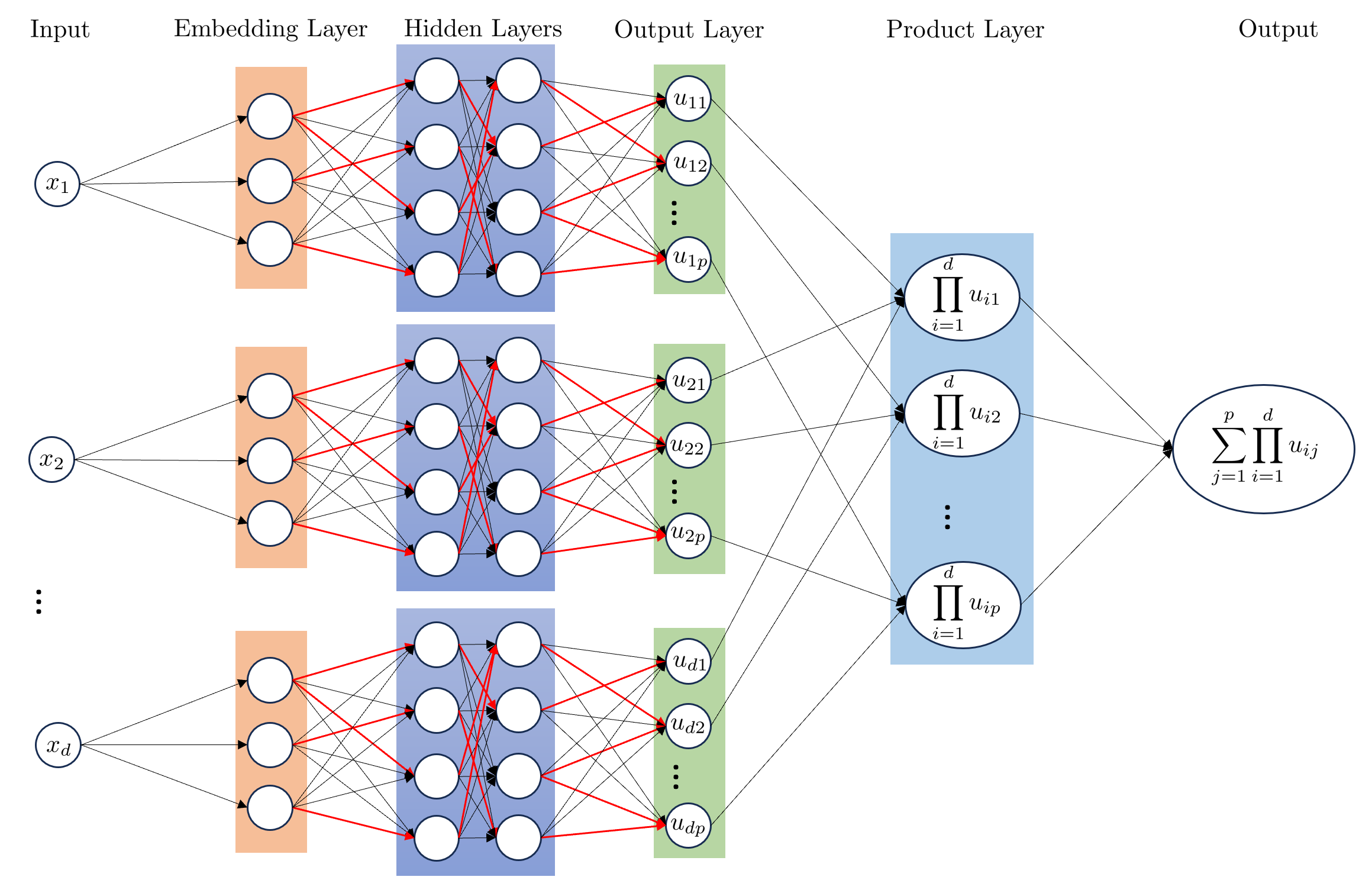}
  \caption{Schematic diagram of the partial evolutionary tensor neural networks (pETNNs). The thick red line represents the updated parameters and other parameters remain unchanged. The embedding layer is used to ensure the boundary conditions.}
  \label{fig2}
\end{figure}

We now provide a detailed introduction to the pETNNs for solving time-dependent PDEs, which is inspired by \cite{PhysRevE.104.045303}. Consider a time-dependent nonlinear partial differential equation of the form,
\begin{equation}\label{3.2.1}
\begin{aligned}
  \frac{\partial u}{\partial t} - \mathcal{N}(u) &= 0, \qquad (x,t) \in \Omega \times [0, T], \\
  u(x, t=0) &= u_0(x)
\end{aligned}
\end{equation}
where $\Omega \subseteq \mathbb{R}^d$ denotes the spatial domain of the computation, $T$ represents the termination time for the simulation, $u=u(x,t)$ is the scalar state function on both space and time and $\mathcal{N}$ is a nonlinear differential operator. Fig. \ref{fig2} shows the schematic diagram of the pETNNs. We add an embedding layer \cite{2020arXiv200610739T} to the original TNN architecture. In this paper, the embedding for periodic boundary conditions is adopted as
$$
x_i \rightarrow a[\cos(b x_i), \sin(b x_i)],
$$
where we use $a$ and $b$ to adjust the computational domain and the periodic frequency, respectively.

It is widely recognized that several prominent techniques for solving time-dependent PDEs with DNNs, such as PINNs \cite{RAISSI2019686}, the deep Galerkin method \cite{SIRIGNANO20181339}, and the deep Ritz method \cite{bing2018DRM}, incorporate both time and spatial variables as inputs and have yielded impressive results. However, the structure of these methods makes it hard to make predictions beyond the training horizon. To forecast solutions without further training and achieve high accuracy, we propose the following pETNNs.

As depicted in Fig.~\ref{fig2} and referred back to \eqref{3.1.1}, the parameters are divided into two parts, time-invariant (thin black lines) and time-variant variables (thick red lines), i.e., $\theta = \{\hat{\theta}(t), \tilde{\theta}\}$, and pETNNs propose to approximate the solution to the PDE \eqref{3.2.1} by
\begin{equation}\label{3.2.2}
 u(x,t) \approx \hat{u}(x, \theta(t))=\sum_{j=1}^p \prod_{i=1}^d \hat{u}_{i, j}\left(x_i , \{\hat{\theta}^i(t), \tilde{\theta}^i\}\right).
\end{equation}
In \eqref{3.2.2}, the set $\{\hat{\theta}^i(t), \tilde{\theta}^i\}$ constitutes the collection of parameters associated with the $i$-th sub-network, where $\hat{\theta}^i(t)$ is designated as a function of time, and $\tilde{\theta}^i$ is fixed over time.

Since we consider the case where we do not have access to the training data to infer the parameter $\theta(t)$ in \eqref{3.2.2} to approximate the PDE solution directly, we propose to develop the governing equation for $\theta(t)$ so that the error is low. Assuming $\hat{\theta}(t)$ is differentiable, the time derivative of solution $\hat{u}$ in \eqref{3.2.2} can be calculated by
\begin{equation*}
  \frac{\partial \hat{u}}{\partial t}=\frac{\partial \hat{u}}{\partial \hat{\theta}}\frac{d \hat{\theta}}{d t}.
\end{equation*}
Inserting the ansatz $\hat{u}$ in the PDE  \eqref{3.2.1}, we follow the Dirac-Frenkel variational principle, e.g.,\cite{PhysRevE.104.045303}, which leads to the following ODE system for time-variant parameters $\hat{\theta}(t)$
\begin{align}\label{3.2.3}
  \frac{d \hat{\theta}}{d t} &= \arg\min_{\gamma} J({\theta}(t), \gamma), \\
  \text{where}\quad J({\theta}(t), \gamma) & = \frac{1}{2}\int_{\Omega}\left|\frac{\partial \hat{u}}{\partial \hat{\theta}}\gamma-\mathcal{N}(\hat{u})\right|^2dx. \nonumber
\end{align}
According to first-order optimality condition, the optimal solution $\gamma_{opt}$ to the minimization problem in \eqref{3.2.3} satisfies the following linear system
\begin{equation}\label{3.2.4}
\left(\int_{\Omega}\left(\frac{\partial \hat{u}}{\partial \hat{\theta}}\right)^{\top}\left(\frac{\partial \hat{u}}{\partial \hat{\theta}}\right) d x\right) \gamma_{o p t}=\int_{\Omega} \left( \frac{\partial \hat{u}}{\partial \hat{\theta}} \right)^{\top} \mathcal{N}(\hat{u}) dx,
\end{equation}
where the entries of the coefficient matrix are given by
\begin{equation}\label{3.2.5}
  \left(\int_{\Omega}\left(\frac{\partial \hat{u}}{\partial \hat{\theta}}\right)^T\left(\frac{\partial \hat{u}}{\partial \hat{\theta}}\right) d x\right)_{i j}=\int_{\Omega} \frac{\partial \hat{u}}{\partial w_i} \frac{\partial \hat{u}}{\partial w_j} d x.
\end{equation}
Here, $w_i, w_{j} \in \hat{\theta}(t)$ are two specific  neural network weights. Followed by the \eqref{3.2.2}, consider all $w \in \hat{\theta}^{k}$ their derivatives in \eqref{3.2.5} can be calculated by
\begin{equation*}
\frac{\partial \hat{u}}{\partial w} 
 =\sum_{\ell_1 =1 }^p \frac{\partial \hat{u}_{k, \ell_{1}}\left( x_k, \{\hat{\theta}^k (t), \tilde{\theta}^k \}\right)}{\partial w} \prod_{\ell_2=1, \ell_2 \neq k}^d \hat{u}_{\ell_2, \ell_1}\left(x_{\ell_2}, \{\hat{\theta}^{\ell_2}(t), \tilde{\theta}^{\ell_2}\}\right),
\end{equation*}
which suggests that the coefficient matrix in \eqref{3.2.4} has a block form. Moreover, similar to \cite{wang2023tensor}, we can efficiently compute the numerical integration involved in the first-order optimality conditions \eqref{3.2.4} due to the variable separation structure of TNN.  We take advantage of the block form of the coefficient matrix and the efficient numerical integration of TNNs to compute the coefficients in \eqref{3.2.4}, which is crucial for high-dimensional problems. We report the details in Appendix~\ref{app:numerical integration}. In practice, the linear system \eqref{3.2.4} is underdetermined, and we use ``linalg.lstsq'' in {\sf NumPy} \cite{harris2020array} to identify $\gamma_{opt}$, which solves \eqref{3.2.4} in the least squares sense with a cut-off ratio for small singular values.

Let $\gamma_{opt}(\hat{\theta})$ denote the least squares solution to \eqref{3.2.4} subject to $\hat{\theta}(t) = \hat{\theta}$, and, in this paper, we consider the predictor-corrector method (also known as the modified-Euler method) to numerically integrate the ODE system \eqref{3.2.3}, which consists of a predictive step and a corrective step.

\noindent $\bullet$ Predictive Step (explicit Euler method):
$$
\hat{\theta}_n^p = \hat{\theta}_n + \Delta t\gamma_{o p t}(\hat{\theta}_n);
$$
\noindent $\bullet$ Corrective Step (implicit Euler method):
$$
\hat{\theta}_{n+1} = \hat{\theta}_n + \frac{\Delta t}{2}\left( \gamma_{o p t}(\hat{\theta}_n)+ \gamma_{o p t}(\hat{\theta}_n^p)\right).
$$
Here $\hat{\theta}_n= \hat{\theta}(n\Delta t)$ with $\Delta t$ the time step size and $n$ the index of time step. In general, the ODE system \eqref{3.2.3} can be numerically integrated using other standard methods.

In the initial phase of the procedure, to generate an initial value for the ODE system \eqref{3.2.3} and the value of time-invariant parameters $\tilde{\theta}$, we train the TNN, $\hat{u}(x, \theta)$, to approximate the initial condition $u(x,t=0) = u_0(x)$. The implementation, including embedded constraints and Dirichlet boundary conditions, follows the methodology outlined in \cite{PhysRevE.104.045303}.
Specifically, we train the tensor neural network with sufficient accuracy by
\begin{equation}\label{eq:TNN_train}
\min_{\theta} \|\hat{u}(x,\theta)-u_0(x)\|.
\end{equation}
Here (so as in Section~\ref{sec:num}) $\|\cdot\|$ represents $L_2$-norm, which can be directly calculated by the conventional Gaussian quadrature schemes in the TNNs architecture. We summarize the entire procedure into Algorithm~\ref{algorithm1}.

\begin{algorithm}[ht]
\caption{pETNNs with Fixed Parameter Update for Solving PDEs}
\label{algorithm1}
\begin{algorithmic}[1]
\Statex \hspace*{-\algorithmicindent} \textbf{Input:} PDE \eqref{3.2.1} with boundary conditions and initial values, Gauss points and quadrature weights, tolerance for the least squares method, time step size $\Delta t$
\Statex \hspace*{-\algorithmicindent} \textbf{Output:} Approximated solutions to the PDE \eqref{3.2.1} represented by the TNNs over $[0,T]$

\State Train the TNNs to approximate $u(x,t=0)=u_{0}(x)$ and obtain $\theta_0$
\State Partition $\theta_0$ into time-dependent $\hat{\theta}_0$ and time-independent $\tilde{\theta}$
\While{$n\Delta t < T$}
    \State Apply Automatic Differentiation (AD) to $\hat{\theta}_n$ and integrate \eqref{3.2.5} using the TNNs structure
    \State Solve for $\gamma_{\text{opt}}(\hat{\theta}_n)$ in \eqref{3.2.4} by the least squares method
    \State Predictive step: $\hat{\theta}_n^p \leftarrow \hat{\theta}_n + \Delta t\gamma_{\text{opt}}(\hat{\theta}_n)$
    \State Apply AD to $\hat{\theta}_n^p$ and integrate \eqref{3.2.5}
    \State Solve for $\gamma_{\text{opt}}(\hat{\theta}_n^p)$ in \eqref{3.2.4} by the least squares method
    \State Corrective step: $\hat{\theta}_{n+1} \leftarrow \hat{\theta}_n + \frac{\Delta t}{2}(\gamma_{\text{opt}}(\hat{\theta}_n) + \gamma_{\text{opt}}(\hat{\theta}_n^p))$
    \State Update and save parameters $\hat{\theta}_{n+1}$
\EndWhile
\State \textbf{return} The updated TNNs with parameters $\{\theta_n\}$
\end{algorithmic}
\end{algorithm}

We want to note that in numerically solving the ODE system \eqref{3.2.3}, the iterative computation unfolds across a series of discrete temporal intervals. Within each interval, we can partition the parameter $\theta$ into time-dependent $\hat{\theta}$ and time-independent $\tilde{\theta}$. Therefore, we present the corresponding Algorithm \ref{algorithm2}. We will examine the performance of these two algorithms and compare them with baseline methods in Section~\ref{sec:compare}. It is noteworthy that Algorithm \ref{algorithm1} may be regarded as a particular instantiation of Algorithm~\ref{algorithm2}, characterized by the constancy of parameters across successive time steps. In general, such a characterization relies on prior knowledge of the PDE system, which may not be available. We close the section with the following remarks.

\begin{algorithm}[ht]
\caption{pETNNs with Randomized Parameter Update for Solving PDEs}
\label{algorithm2}
\begin{algorithmic}[1]
\Statex \hspace*{-\algorithmicindent} \textbf{Input:} A PDE \eqref{3.2.1} with boundary conditions and initial values, Gauss points and quadrature weights, tolerance for the least squares method, time step size $\Delta t$,
\Statex \hspace*{-\algorithmicindent} \textbf{Output:} Approximated solutions to the PDE \eqref{3.2.1} represented by the TNNs over $[0,T]$

\State Train the TNNs to approximate $u(x,t=0)=u_{0}(x)$ and obtain $\theta_0$
\While{$n\Delta t < T$}
    \State Partition $\theta_n$ randomly into time-dependent $\hat{\theta}_n$ and time-independent $\tilde{\theta}_n$
    \State Apply AD to $\hat{\theta}_n$ and integrate \eqref{3.2.5} using the TNNs structure
    \State Solve for $\gamma_{\text{opt}}(\hat{\theta}_n)$ in \eqref{3.2.4} by the least squares method
    \State Predictive step: $\hat{\theta}_n^p \leftarrow \hat{\theta}_n + \Delta t\gamma_{\text{opt}}(\hat{\theta}_n)$
    \State Apply AD to $\hat{\theta}_n^p$ and integrate \eqref{3.2.5}
    \State Solve for $\gamma_{\text{opt}}(\hat{\theta}_n^p)$ in \eqref{3.2.4} by the least squares method
    \State Corrective step: $\hat{\theta}_{n+1} \leftarrow \hat{\theta}_n + \frac{\Delta t}{2}(\gamma_{\text{opt}}(\hat{\theta}_n) + \gamma_{\text{opt}}(\hat{\theta}_n^p))$
    \State Update and save parameters $\theta_{n+1}=\{\hat{\theta}_{n+1},\tilde{\theta}_n\}$
\EndWhile
\State \textbf{return} The updated TNNs with parameters $\{\theta_n\}$
\end{algorithmic}
\end{algorithm}

\begin{remark} \label{rem1}
  It is widely acknowledged that DNNs exhibit inherent parameter redundancy.
  The redundancy not only bolsters robustness but also enables us to attain results that match or surpass those of full updates through updating partial parameters, either randomly or with careful selection. However, considering that different parameters indeed play distinct roles within the neural network, manually specifying sampling criteria can still provide us with additional benefits, such as enhancing the algorithm's stability over extended periods. {  In general, the optimal partition of parameter into time-dependent and time-independent is hard to achieve and case-by-case. We test several strategies in Section~\ref{sec:num}.}
\end{remark}

{ 
\subsection{Relations to the existing methods}

In Section~\ref{sec:pETNN}, we have proposed the pETNNs for solving time-dependent PDEs and introduced the corresponding algorithms. It is essential to clarify how our approach differs from existing methods and highlight the specific improvements achieved by pETNNs.

\noindent \textbf{The TNN-based methods for PDEs} The TNN architecture used in pETNNs is identical to a number of existing TNN-based methods for PDEs, e.g., \cite{wang2023tensor}. The evolutionary parametric approximation distinguishes pETNNs from others by using ODE systems derived from the variational principle to drive the parameters, rather than discretizing the spatial-temporal domain and minimizing specific loss functions for each mesh. Leveraging the representation capabilities of TNN, pETNNs also offer physical interpretations for the parameters, helping to prevent potential training failures.

\noindent \textbf{Hierarchical tensor methods for PDEs} The use of TNNs in our parametrization links pETNNs to the hierarchical tensor method employed for solving high-dimensional PDEs, as evidenced in \cite{Markus2016Tensor, bachmayr2023low}. This approach addresses the notorious challenge of the curse of dimensionality encountered in the numerical resolution of high-dimensional PDEs. Inspired by the classical technique of variable separation and tensor decompositions, hierarchical tensor methods strive for low-rank approximations of the solution tensor of the PDE system, depending largely on the choices of spatial discretization and tensor decomposition format. Further details can be found in the review papers \cite{Markus2016Tensor, falco2019dirac}. Distinctly, the TNNs described in \eqref{3.1.1} offer a continuous representation of general rank-$p$ tensors, not limited to any specific tensor decomposition format. As a consequence, pETNNs' implementation is independent of tensor decomposition algorithms, enhancing its storage efficiency and facilitating compatibility with high-dimensional PDEs with minimal prior knowledge of the solution structure. Moreover, in pETNNs, low-rank structures are integrated not only into the network architecture of the TNNs parameterization but also into the updating strategy. 

\noindent \textbf{The evolutionary deep neural network} In the field of AI for PDEs, there exist a variety of methodologies that make use of evolutionary parametric approximation techniques as in pETNNs. We want to clarify that using the Dirac-Frenkel variational principle to develop nonlinear evolutionary parametrization is a conventional method in molecular dynamics and computational chemistry. Recent works e.g., \cite{BRUNA2024112588}, have leveraged the variational principle to derive systems of nonlinear evolution equations for the parameters of the parametric representation like DNN. We refer the readers to \cite{BRUNA2024112588, PhysRevE.104.045303} and references therein for the details. We adapt these frameworks to TNNs for solving high-dimensional time-dependent PDEs. The main innovation of pETNNs lies in their low-rank parameterization and numerical implementations, such as the integration of the coefficient matrix, enabling them to tackle high-dimensional challenges.

The efficacy of this approach will be demonstrated through numerical tests in Section~\ref{sec:num}. At a more advanced level, our method can be classified within the realm of low-rank approximation strategies for dynamical systems.
}

\subsection{A Posterior Error Estimation of pETNNs}

In this section, we discuss a posterior error bound of pETNNs subject to spatial discretization. Although pETNNs belong to the class of regularized dynamical parametric approximation \cite{feischl2024regularized}, we cannot directly adopt the error estimation result in \cite{feischl2024regularized} since $\mathcal{N}(u)$ in \eqref{3.2.1}, in general, does not satisfy the Lipschitz-continuous condition. To resolve the issue, we will further assume $\mathcal{N}(u)$ is linear and \eqref{3.2.1} can be transformed into a system of tensor-valued ODEs
\begin{equation}\label{eq:ODE}
    \frac{d \bm{u}}{d t} = \bm{N} \bm{u}, \quad \bm{u}(0) = \bm{u}_0,
\end{equation}
where $\bm{u}: [0, T] \rightarrow \mathbb{R}^{n_1 \times n_2 \times \dots \times n_d}$ is a multi-dimensional array of real numbers (the solution tensor), $\bm{N}$ is a finite-dimensional linear operator (the discrete form of $\mathcal{N}$). In general, the structure of $\bm{N}$ depends on the spatial discretization of $u$ in \eqref{3.2.1}, as well as on the tensor format utilized for $\bm{u}$ \cite{Daniele1}. For example, when the spatial domain is a $d$-dimensional box, the solution tensor $\bm{u}$ can be defined by approximating $u(x,t)$ using a finite-dimensional Galerkin basis with $N$ degrees of freedom in each direction yields a total number of $N^{d}$ degrees of freedom, and $\bm{u}_0$ in \eqref{eq:ODE} is the projection of the initial condition $u_{0}$ in \eqref{3.2.1}.

Another motivation for considering spatial discretization comes from the numerical implementations of pETNNs. Notice that we propose to solve the governing equation of $\hat{\theta}(t)$ \eqref{3.2.3} via the linear system \eqref{3.2.4} whose coefficients are determined by direct numerical integrations of $\hat{u}$ and its partial derivatives on $\Omega$ (See Appendix~\ref{app:numerical integration} for the details). Therefore, even though the formulation of pETNNs is mesh-free and does not rely on a specific spatial discretization, numerically, the dynamics of $\hat{\theta}(t)$ governed by \eqref{3.2.3} is determined by the information $\hat{u}(x, \theta(t))$ over grid points.

Subject to the same spatial discretization in \eqref{eq:ODE}, let $\hat{\bm{u}}$ be the solution tensor of the approximated solution $\hat{u}(x, {\theta}(t))$ generated by pETNNs. Then, $\hat{\bm{u}}$ is governed by
\begin{equation}\label{eq:ODE_approx}
    \frac{d \hat{\bm{u}}}{d t} = \bm{N} \hat{\bm{u}} + \bm{r}, \quad \hat{\bm{u}}(0) = \hat{\bm{u}}_0,
\end{equation}
where $\hat{\bm{u}}_0$ is the projection of the pETNNs initialization $\hat{u}(x, \theta_{0})$ obtained from the TNN training in \eqref{eq:TNN_train} and $\bm{r}$ denotes the defect introduced by pETNNs. We assume the defect satisfies $\|\bm{r}(t)\|_{D} \leq \delta(t)$, where $\|\cdot\|_{D}$ is some norm corresponding to the spatial discretization. Consequently, we reduce the problem to a perturbation analysis between two linear ODE systems \eqref{eq:ODE} and \eqref{eq:ODE_approx}. While $\delta(t)$ can be monitored during the computation and is thus available in a posterior manner. We propose the following posterior error estimation of pETNNs.

\begin{proposition}\label{prop:error_bound}
Consider the linear ODE system \eqref{eq:ODE} and its approximation \eqref{eq:ODE_approx} satisfying
\begin{enumerate}
    \item (bounded linear operator) $\|\bm{N}(\bm{u} -\bm{v})\|_{D} \leq L \|\bm{u} -\bm{v} \|_{D}$, and
    \item (defect bound) $\|\bm{d}(t)\|_{D} \leq \delta(t)$, 
\end{enumerate}
We have the error bound
\begin{equation*}
\|\bm{u}(t)-\hat{\bm{u}}(t)\|_{D} \leq e^{L t} \|\bm{u}_0 - \hat{\bm{u}}_0\|_{D} + \int_{0}^{t}e^{L (t-s)}\delta (s) d s.
\end{equation*}
\end{proposition}

\begin{proof}
At $t=0$, we have $\| \bm{u}(0) - \hat{\bm{u}}(0) \|_{D} = \|\bm{u}_0 - \hat{\bm{u}}_0\|_{D}$, and the perturbation introduced by pETNNs is bounded by
\begin{equation*}
\left\| \frac{d}{d t} (\bm{u} - \hat{\bm{u}}) \right\|_{D} \leq \| N(\bm{u} - \hat{\bm{u}})\|_{D} + \|\bm{r}\|_{D} \leq L \| \bm{u} - \hat{\bm{u}}\|_{D} + \delta(t).
\end{equation*}
The standard Growwall's inequality yields the error bound.
\end{proof}
We want to remark on the assumptions and constants in Prop.~\ref{prop:error_bound}. To begin with, we can take the Lipschitz constant $L$ as the spectral radius of $\bm{N}$. The relation between $L$ and the dimension $d$ is known as the stiffness in high-dimensional PDEs. For further discussion, we refer the reader to \cite{Daniele1} and references therein. In \cite{feischl2024regularized}, the defect bound is characterized by the regularized least squares loss of the linear system \eqref{3.2.4}. Here, we consider a different regularization approach by truncating the small singular values of the coefficient matrix in \eqref{3.2.4}, and introduce the defect bound posteriorly. Lastly, notice that $\bm{u}_{0}$ and $\hat{\bm{u}}_{0}$ are projections of the initial condition $u_{0}$ and $\hat{u}(x,\theta_{0})$, respectively, and we have $\|\bm{u}_0 - \hat{\bm{u}}_0\|_{D}  \lesssim \|\hat{u}(x,\theta_0)-u_0(x)\| $, which is controlled by the approximation error of TNNs thanks to the universal approximation theory of TNNs \cite{wang2023tensor}.





\section{Numerical Experiments}
\label{sec:num}

We present numerous numerical examples designed to demonstrate the efficacy of the pETNNs. We commence by examining the performance of various parameter update strategies within the pETNNs, including whole parameters update, and partial parameters update (Algorithm~\ref{algorithm1} and Algorithm~\ref{algorithm2}). To illustrate the superiority of the pETNNs in accuracy, comparative experiments delineating their performance relative to the EDNNs are conducted. The long-term stability of the pETNNs is also demonstrated by an experiment with an ultra-long time. Additionally, we intend to employ the pETNNs for the approximate solutions of time-dependent complex PDEs with periodic and homogeneous Dirichlet boundary conditions, such as the 2D incompressible Navier-Stokes equations, the 10D heat equation, the 10D transport equation, and the 10D and 20D simplified Korteweg-de Vries (KdV) equation. {  All experiments were performed on an NVIDIA A100 GPU (80GB) and an Intel Xeon Gold 6326 CPU (2.90GHz).}

\subsection{Comparisons and parameter update strategies}\label{sec:compare}

In the first experiment, we use the following transport equation
$$
\frac{\partial u}{\partial t} + \sum_{i=1}^d \frac{\partial u}{\partial x_i}=0, \quad x \in \Omega=[-1,1]^d
$$
to compare our method with baselines, including the finite difference methods, PINNs, EDNNs and ETNNs. We also test Algorithm~\ref{algorithm1} and Algorithm~\ref{algorithm2} with partial and full updates to demonstrate the parameter update strategies. The considered transport equation reads as with the initial values
$$
u(x, t=0)=\prod_{i=1}^d \sin \left(\pi x_i\right)
$$
and periodic boundary conditions. The analytical solutions are
$$
u(x, t)=\prod_{i=1}^d \sin \left(\pi\left(x_i- t\right)\right),
$$
where $d=2,3.$ The absolute error $\varepsilon_{a}$ and relative error $\varepsilon_{r}$ are respectively defined as
$$
\varepsilon_{a}=\|\hat{u}(t)-u(t)\|, \qquad \varepsilon_{r}=\frac{\|\hat{u}(t)-u(t)\|}{\|u(t)\|}.
$$

As for the network architecture, we fix the hyperbolic tangent, $\mathrm{tanh}$, as the activation function for all models. For the PINNs in the 2D case, we use a 6-layer fully-connected neural network with 20 neurons per layer. To impose the periodic boundary condition, the network input $x_j$ is replaced with $\{\cos(\pi x_j), \sin(\pi x_j)\}$. We use a time-marching strategy\cite{wight2020solving} to split the time domain into intervals of length $2$, and in each spatial-temporal domain of the form $[-1,1]^2\times [t_0, t_0+2]$, we create a uniform mesh of size $32^2\times 400$. The PINNs model is trained via full-batch gradient descent using the Adam optimizer\cite{kingma2014adam} for $10^6$ iterations and followed by L-BFGS \cite{byrd1995limited} for $5\times 10^4$ iterations.

For the EDNNs, the architecture in the 2D case includes three hidden layers with 20 neurons per layer, whereas in the 3D case, it has four hidden layers with 32 neurons per layer. In contrast, the ETNNs maintain a consistent structure across 2D and 3D cases, with each sub-network containing two hidden layers of 20 neurons each and a fixed time step $1e-3$.

{  In the 2D case, we also consider the finite difference method as a baseline. In the spatial domain, we implement the fifth-order weighted essentially non-oscillatory (WENO) scheme\cite{Titarev2004, Hesthaven2018} with a uniform mesh of size $480\times 480$, corresponding to the 480 Gaussian quadrature points used in the numerical integration in pETNNs. In the time evolution, we use the predictor-corrector (modified-Euler) method with a time step $1e-3$. Thus, the CFL number is $0.24$. For detailed discussions on WENO scheme and finite difference methods for transport equations, we refer to \cite{Shu1998, GREC2025275}.
}

Table~\ref{table1} shows the relative errors of the finite difference method, PINNS, pEDNNs (EDNNs with partial update strategy, referred to as the randomized sparse Neural Galerkin schemes \cite{NEURIPS2023_0cb310ed}) and pETNNs for solving the 2D transport equation, and Table \ref{table3} enumerates the relative errors for their 3D counterparts. As suggested by the error statistics, pETNNs have a relatively higher accuracy due to the special structure of the tensor neural networks. {  The accuracy of pETNNs suggests that the evolutionary TNNs \eqref{3.2.2}, under the selected architecture, is a reasonable hypothesis space for 2D transport equations, which inspires us to apply time-marching strategy to TNNs. By treating the time variable as the $(d+1)$th input, we formulate
the space-time TNNs and train in the same way as the PINNs. The performance is considerably lower than that of pETNNs, with numerical results for space-time TNNs provided in Appendix~\ref{app:st_TNN}. 
}
\begin{table}[htbp]
\centering
\begin{tabular}{@{}ccccccc@{}}
\toprule
& \multicolumn{6}{c}{relative error $\varepsilon_{r}$}                                             \\ \cmidrule(l){2-7}
     & t=0       & t=2       & t=4       & t=6       & t=8       & t=10      \\\cmidrule(l){2-7}
{  WENO5} &  0  & 5.85E-05  & 1.17E-04  & 1.75E-04 & 2.34E-04 & 2.92E-04 \\
PINNs &  3.67E-05 & 2.10e-04 & 3.28E-04 & 3.14E-04 & 3.89E-04 &  6.87E-04\\
EDNNs(1/3) & 1.14E-04 & 1.15E-04 & 1.18E-04 & 1.23E-04 & 1.30E-04 & 1.39E-04\\
EDNNs(1/2) & 1.14E-04 & 1.15E-04 & 1.16E-04 & 1.20E-04 & 1.30E-04 & 1.25E-04\\
EDNNs(all) & 1.14E-04 & 1.15E-04 & 1.16E-04 & 1.17E-04 & 1.20E-04 & 1.23E-04 \\ 
ETNNs(1/3) & 1.79E-05 & 1.91E-05 & 2.20E-05 & 2.64E-05 & 3.07E-05 & 3.83E-05\\
ETNNs(1/2) & 1.79E-05 & 1.89E-05 & 2.00E-05 & 2.13E-05 & 2.36E-05 & 2.66E-05\\
ETNNs(all) & 1.79E-05 & 1.96E-05 & 2.36E-05 & 2.66E-05 & 3.24E-05 & 3.56E-05 \\ \bottomrule 
\end{tabular}
\caption{Relative errors of the fifth order WENO scheme finite difference method (WENO5), PINNs, pEDNNs and pETNNs with different parameter update strategies for solving the 2D transport equation, ($*$) denotes that only $*$ of parameters are updated. }
\label{table1}
\end{table}



We also compare the results of different parameter update strategies. It shows that the efficacy of randomly updating a subset of TNN parameters is comparable to that of updating the entire set of parameters, and even randomly updating $1/6$ of parameters can achieve excellent results, including numerical accuracy and stability. In contrast, as for the performance of pEDNNs, even though we reported the best relative error we obtained over various hyperparameter values in Table~\ref{table3}, the partial update strategy is still suboptimal.

\begin{table}[ht]
\centering
\begin{tabular}{@{}ccccccc@{}}
\toprule
& \multicolumn{6}{c}{relative error $\varepsilon_r$}                                     \\ \cmidrule(l){2-7}
& t=0       & t=1       & t=2       & t=3       & t=4       & t=5      \\ \cmidrule(l){2-7}
EDNNs(1/2) & 2.87E-04 & 9.11E-03 & 1.22E-02 & 1.72E-02 & 2.03E-02 & 2.32E-02\\
EDNNs(all) & 2.87E-04 & 2.88E-04 & 2.87E-04 & 2.89E-04 & 2.89E-04 & 2.91E-04  \\
ETNNs(1/6) & 3.23E-05 & 4.53E-05 & 5.87E-05 & 7.52E-05 & 8.95E-05 & 1.04E-04 \\
ETNNs(1/3) & 3.23E-05 & 3.31E-05 & 3.47E-05 & 3.70E-05 & 3.85E-05 & 4.02E-05 \\
ETNNs(1/2) & 3.23E-05 & 3.26E-05 & 3.33E-05 & 3.42E-05 & 3.54E-05 & 3.73E-05 \\
ETNNs(2/3) & 3.23E-05 & 3.26E-05 & 3.30E-05 & 3.38E-05 & 3.45E-05 & 3.57E-05 \\
ETNNs(5/6) & 3.23E-05 & 3.26E-05 & 3.31E-05 & 3.39E-05 & 3.44E-05 & 3.56E-05 \\
ETNNs(all) & 3.23E-05 & 3.27E-05 & 3.32E-05 & 3.38E-05 & 3.49E-05 & 3.61E-05 \\ \bottomrule
\end{tabular}
\caption{Relative errors of pEDNNs and pETNNs with different parameter update strategies for solving the 3D transport equation, ($*$) denotes that only $*$ of parameters are updated.}
\label{table3}
\end{table}

\begin{figure}[ht!]
  \centering
  \includegraphics[width=0.45\textwidth]{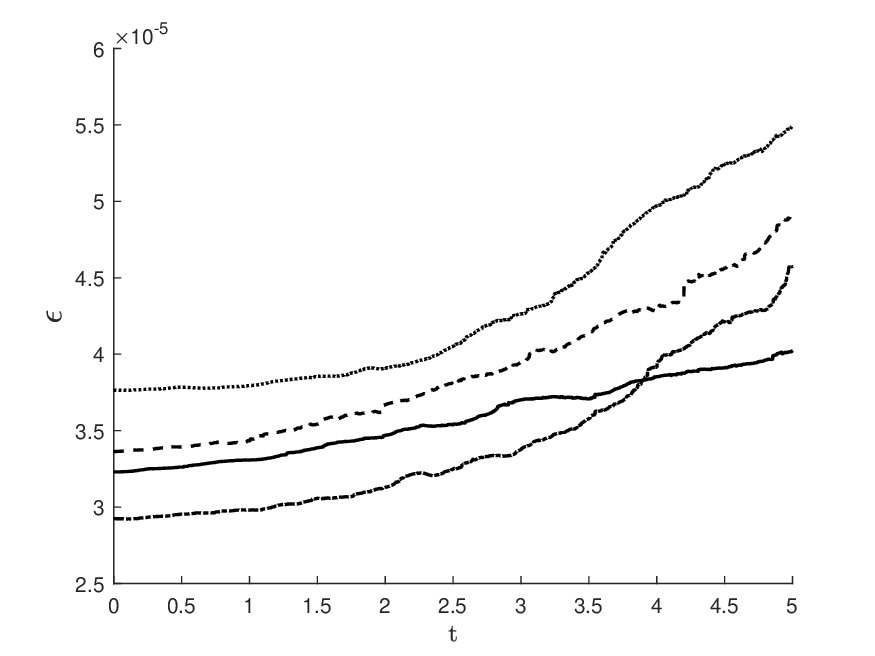}
  \caption{The averaged absolute errors of 3D transport equations over independent runs of Algorithm~\ref{algorithm2}. We updated 1/3 of the parameters randomly per iteration. The four error curves correspond to four different initial $\theta_0$ configurations.}\label{Exp3.1a}
\end{figure}

We further test the robustness of Algorithm~\ref{algorithm2} under different initial parameters $\theta_0$ obtained by independent training procedures of TNNs to approximate the initial condition. Starting with different satisfactory initial parameters, $\theta_0$, we update 1/3 of the parameters randomly at each iteration. The averaged absolute errors, calculated as the mean of results from three independent runs of Algorithm~\ref{algorithm2}, (to reduce the uncertainty introduced by randomly selecting parameters) are illustrated in Fig.~\ref{Exp3.1a}. The narrow error band (about $2.5\times 10^{-5}$) indicates a high degree of methodological stability with respect to the initial parameters.

\begin{figure}[ht]
  \centering
  \includegraphics[width=0.45\textwidth]{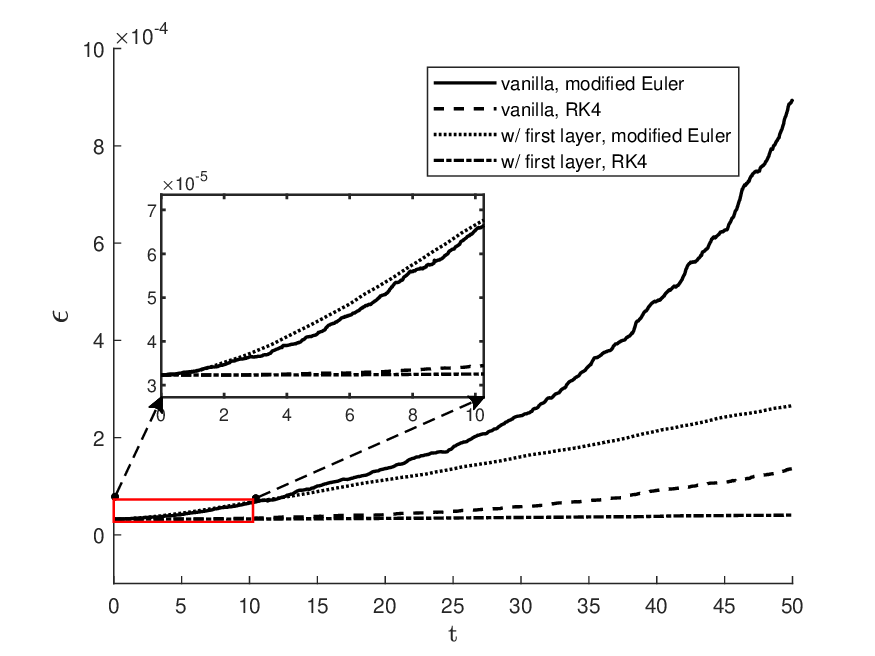}
  \includegraphics[width=0.45\textwidth]{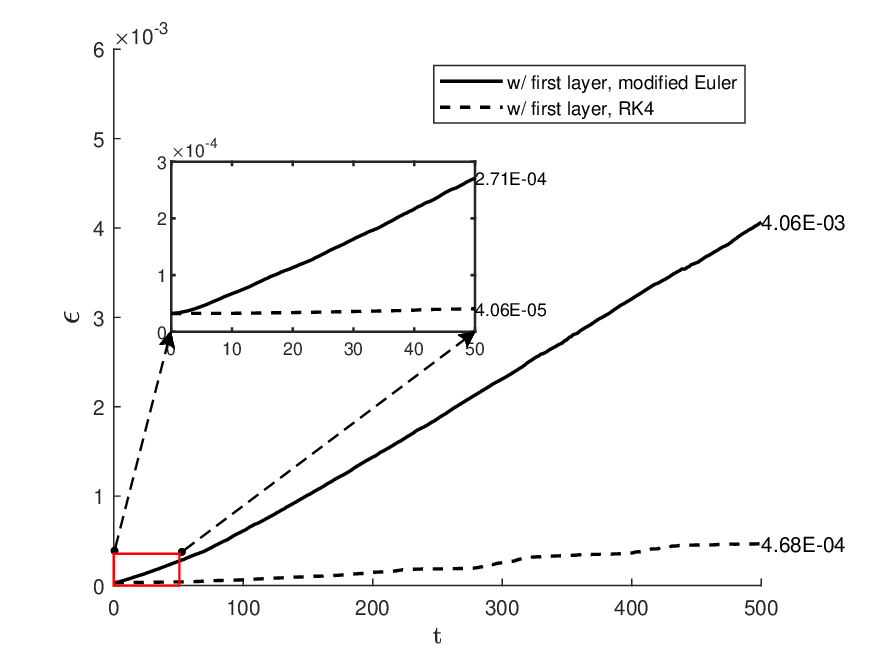}
  \caption{Averaged absolute errors of updating 1/3 of the parameters randomly per iteration: 50 (left) and 500 (right) time units. $``vanilla"$ represents random update in every time step; $``w/~first~layer"$ represents that the inclusion of the update parameters in the first layer is mandatory; $``modified~ Euler"$ represents predictor-corrector (modified-Euler) method, and $``RK4"$ represents explicit fourth-order Runge-Kutta method.}
  \label{Exp3.1b}
\end{figure}

To assess the long-time stability of Algorithm~\ref{algorithm2} and the cumulative numerical errors over extended periods, we extend the time interval to $[0,50]$, and consider the fourth order Runge-Kutta method (RK4) as the control group for comparison in our study. We still update $1/3$ of the parameters per iteration, where the averaged absolute error is calculated as mean values from three independent runs. Our parameter update strategy employs two distinct approaches: the first consistently incorporates the first layer in update (denoted by ``w/ first layer"), which involves all the weights that connect the embedding layer and hidden layer in Fig.~\ref{fig2}, while the second imposes no mandatory inclusion of any specific layer (denoted by ``vanilla"). We draw the results in Fig.~\ref{Exp3.1b}(left), and even by time $t=50$, the averaged absolute error of RK4 remains on the order of $10^{-5}$, demonstrating the method's sustained stability over a long time. Furthermore, recognizing the potential superiority of the ``with first layer" strategy for parameter updates, we replicate the experimental framework to conduct an additional experiment over an ultra-long time, $[0,500]$. The results of this experiment, which are depicted in Fig.~\ref{Exp3.1b}(right), underscore our algorithm's exceptional accuracy and remarkable stability over an extensively prolonged period.

\begin{figure}[ht]
  \centering
  \includegraphics[width=0.45\textwidth]{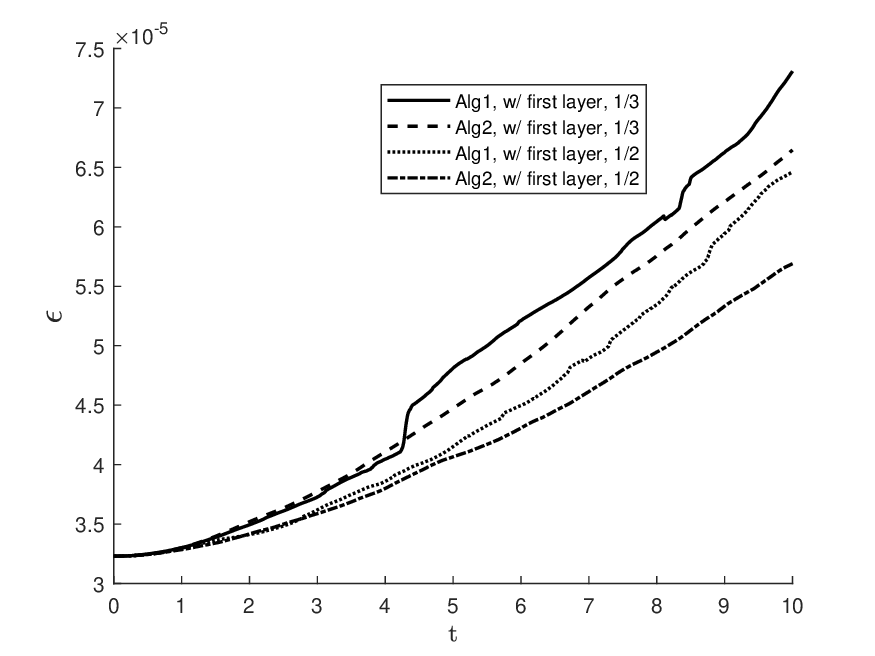}
  \caption{Averaged absolute errors of different update strategies. $``w/~first~layer"$ represents that the inclusion of the update parameters in the first layer is mandatory, and $1/3(1/2)$ represents that only  $1/3(1/2)$  of parameters are updated at each time step.}\label{Exp3.1d}
\end{figure}

The main difference between Algorithm \ref{algorithm1} and Algorithm \ref{algorithm2} is that Algorithm \ref{algorithm2} resamples the parameters to be evolved at each time step. In Fig. \ref{Exp3.1d}, we plot the absolute errors of both algorithms with different parameter update strategies. We find that resampling the parameters to be evolved at each step always outperforms constantly updating a fixed subset of parameters. Thus, we adopt Algorithm~\ref{algorithm2} in the following numerical experiments.

\subsection{The Navier-Stokes Equations} \label{sec:NS}

The incompressible Navier-Stokes (NS) equations are a set of nonlinear PDEs that describe the motion of fluid substances such as liquids and gases. Following \cite{li2022learning}, we consider the two-dimensional Kolmogorov flow for a viscous, incompressible fluid with periodic boundary condition
\begin{equation}\label{eq:NS}
\frac{\partial \bm{q}}{\partial t} + (\bm{q} \cdot \nabla) \bm{q} = -\nabla p + \frac{1}{Re} \nabla^2 \bm{q} + \bm{f}, \quad  \nabla \cdot \bm{q} = 0, \quad (x,y) \in [0, 2\pi]^2,  
\end{equation}
where $\bm{q} = (u(x,y.t),v(x,y,t))$ denotes the velocity, $p$ is the pressure, $Re>0$ is the Reynolds number, and $\bm{f} = (\sin(4y), 0)^{\top}$ is the source term. Unlike the other equations tested in the paper, e.g., transport and heat equations, the explicit form of the exact solution is unavailable. In \cite{li2022learning}, the exact solution data\footnote{The dataset is available at https://zenodo.org/records/7495555.} is obtained by solving the equation in vorticity form using the pseudo-spectral split-step method \cite{chandler2013invariant}, and we adopt their exact solution data for $Re=40$ in our test.

\begin{figure}[ht!]
\centering
\subfigure[pETNNs, $t=0$, $p=20$]{ \includegraphics[width=.225\linewidth]{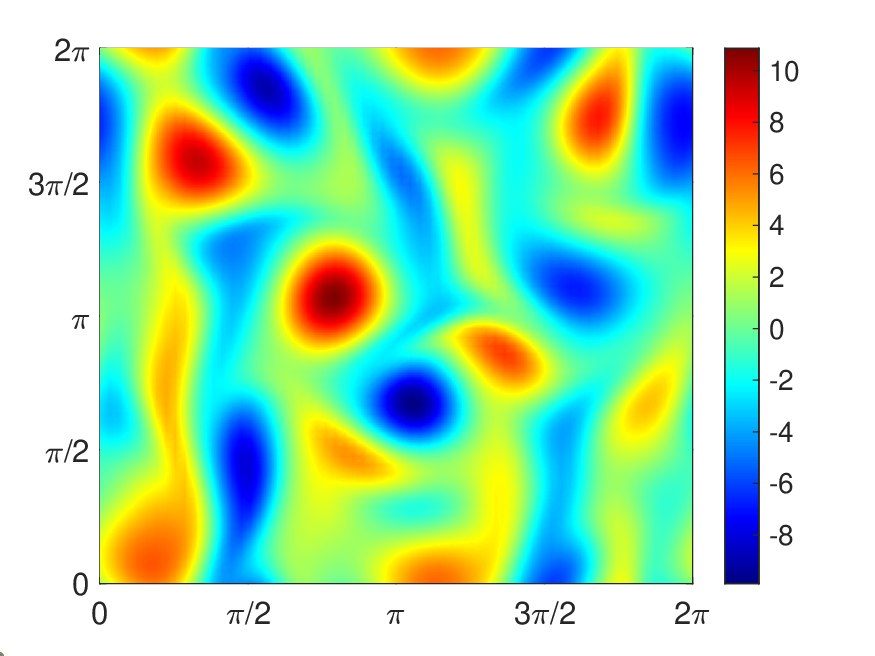}}
\subfigure[Error, $t=0$, $p=20$]{ \includegraphics[width=.225\linewidth]{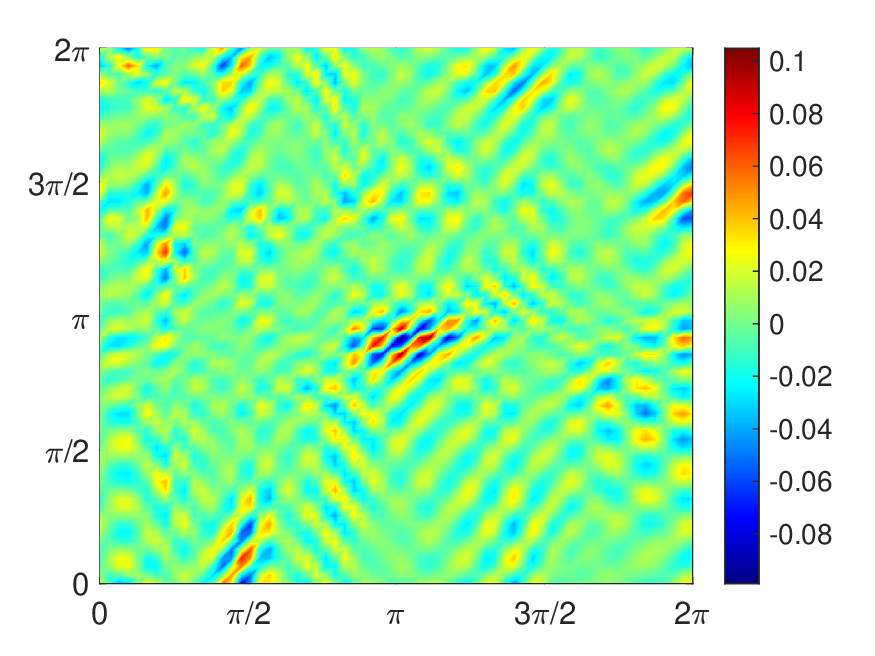}}
\subfigure[pETNNs, $t=10$, $p=20$]{ \includegraphics[width=.225\linewidth]{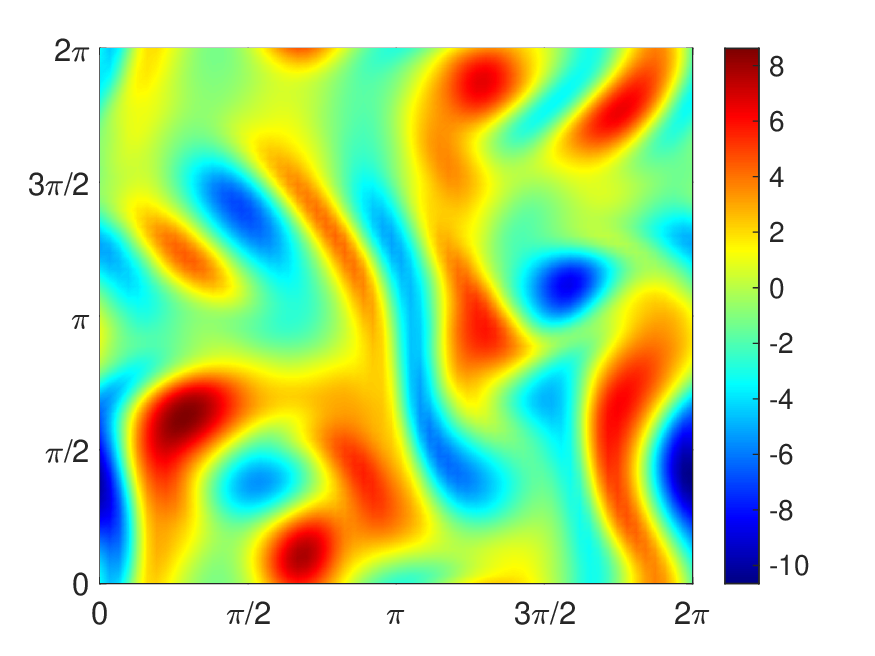}}
\subfigure[Error, $t=10$, $p=20$]{ \includegraphics[width=.225\linewidth]{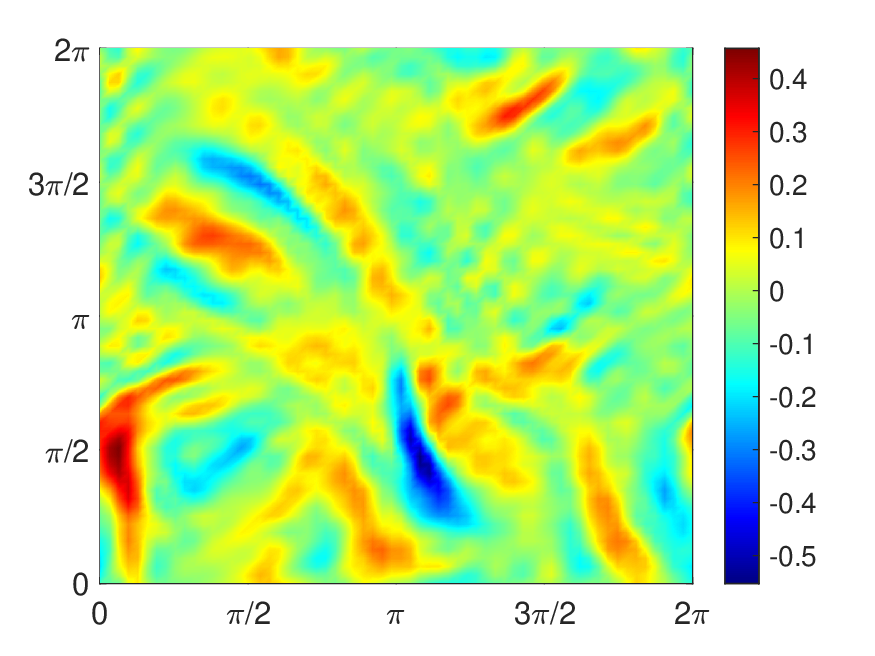}}

\subfigure[pETNNs, $t=0$, $p=30$]{ \includegraphics[width=.225\linewidth]{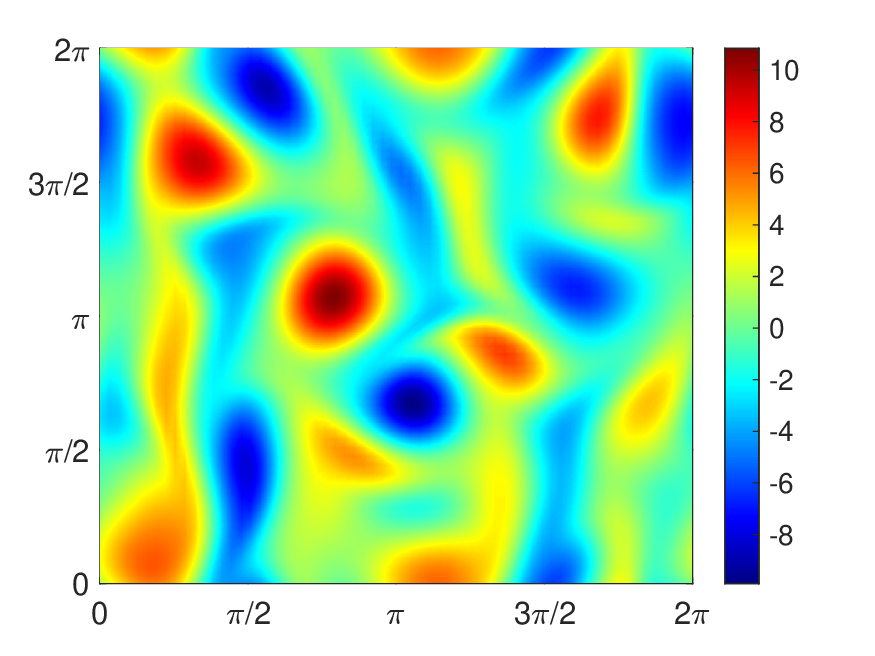}}
\subfigure[Error, $t=0$, $p=30$]{ \includegraphics[width=.225\linewidth]{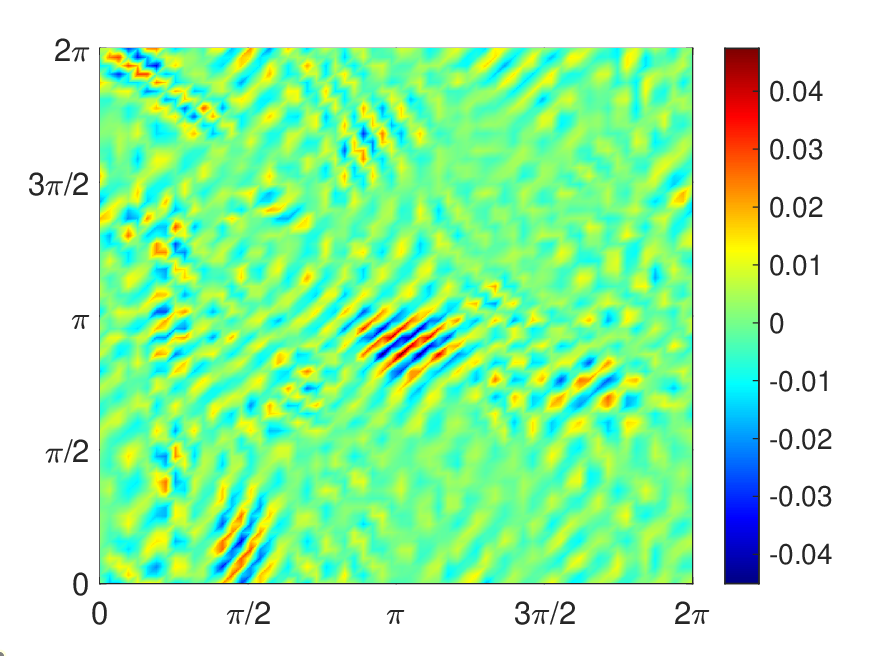}}
\subfigure[pETNNs, $t=10$, $p=30$]{ \includegraphics[width=.225\linewidth]{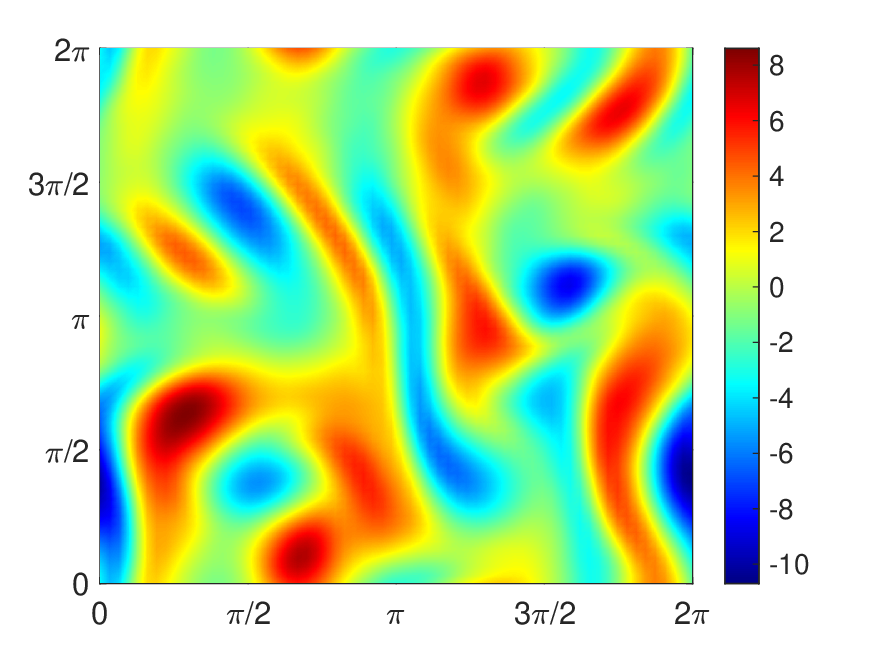}}
\subfigure[Error, $t=10$, $p=30$]{ \includegraphics[width=.225\linewidth]{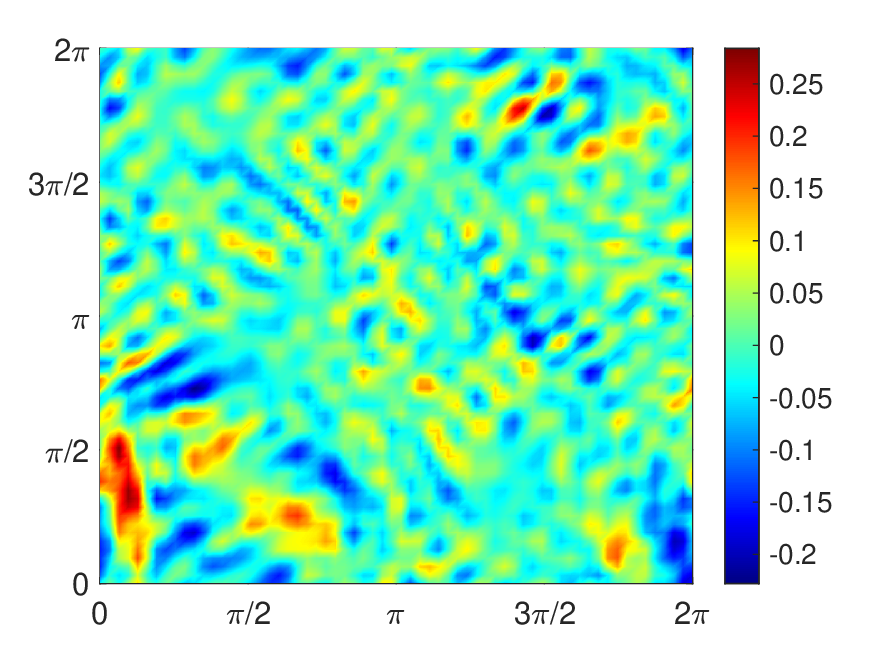}}
\caption{The pETNNs solution and error profiles to the Navier-Stokes equations  and the error at time $t=0$ and $t=10$. Figs.~(a)-(d) correspond to the pETNN solution of rank $p=20$, while Figs.~(e)-(h) correspond to the pETNN solution of rank $p=30$. The error profile is defined as the difference between the pETNNs solution $(\hat\psi_{y}, -\hat\psi_{x})$ and the reference solution generated by the Fourier-Galerkin approach.}
\label{Exp3.2a}
\end{figure}

To preserve the divergence-free condition in the NS equation, we assume that $u=\psi_y, v=-\psi_x$ for some function $\psi(x,y,t)$, and we are trying to solve for $\psi(x,y,t)$ in this experiment. In particular, by substituting $(u,v)$ in \eqref{eq:NS} with $(\psi_{y}, -\psi_{x})$, the PDE \eqref{eq:NS} transforms into a PDE system for $\psi$ of the form
\begin{equation}\label{eq:psi}
\frac{\partial}{\partial t} \psi_{y} = \mathcal{N}_{1}(\psi_{x}, \psi_{y}), \quad \frac{\partial}{\partial t} \psi_{x} = \mathcal{N}_{2}(\psi_{x}, \psi_{y}),
\end{equation}
where 
\begin{equation*}
\begin{split}
\mathcal{N}_{1}(\psi_{x}, \psi_{y}) & = -p_{x} - (\psi_{y}^2)_{x} + (\psi_{x}\psi_{y})_{y} + \frac{1}{Re} \Delta \psi_y, \text{ and} \\
\mathcal{N}_{2}(\psi_{x}, \psi_{y}) & = p_{y} + (\psi_{x}^2)_{y} - (\psi_{x}\psi_{y})_{x} + \frac{1}{Re} \Delta \psi_x.
   \end{split}
\end{equation*}
Introduce a TNN parameterization to $\psi$, $\psi \approx \hat{\psi}(x,y,\theta(t))$, and we apply pETNNs to the two time-dependent PDEs in \eqref{eq:psi} and couple the resulting ODE systems of $\theta(t)$. Such auxiliary function $\psi$ can be derived by the vorticity of the flow, $\omega = v_{x} - u_{y}$, which satisfies $\omega = - \Delta \psi$, that is, $\psi$ is a solution to the Poisson equation with the vorticity as the source term. We use the relation to generate the initial conditions $\hat{\psi}(x,y,\theta_{0})$ for pETNNs. Instead of training the TNN to fit an explicit initial condition like \eqref{eq:TNN_train}, we identify $\theta_{0}$ by solving the Poisson equation based on the TNN parameterization of $\psi$. As for the vorticity data, the dataset in \cite{li2022learning} consists of realizations of vorticity trajectories on the time interval $[0, 500]$ based on random initialization. We use the vorticity information at $t=10$ to formulate the Poisson equation, and solve the PDE by PINNs via energy natural gradient descent \cite{Muller2023achieving}.

In this experiment, for the two sub-networks in the TNNs, we adopt $\mathrm{tanh}$ activation and set two hidden layers with 30 neurons in each layer. We apply Algorithm~\ref{algorithm2} with RK4 scheme and update $1/2$ of the parameters in each iteration with a fixed time step $5e-3$. {  Fig.~\ref{Exp3.2a} reports the solution and error profiles at $t=0$ and $t=10$ of pETNNs for two rank setups $p=20$ and $p=30$. The runtime of pETNNs to reach the 10 time units are 14,730s and 58,758s for $p=20$ and $p=30$, respectively.}

\begin{figure}[ht!]
  \centering
  \includegraphics[width=0.45\textwidth]{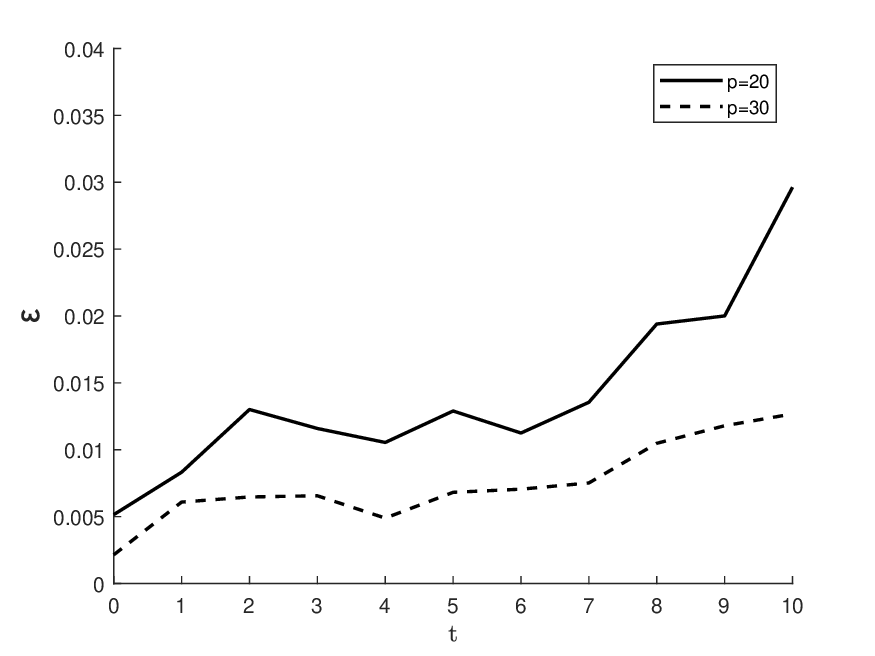}
  \caption{The relative $L^{2}$-error of the vorticity $\omega$ up to $10$ time unites: pETNN with rank $p=20$ (solid) and $p=30$ (dash). The step size is one time unit due to the exact solution data from \cite{li2022learning}.}
  \label{Exp3.2b}
\end{figure}

The solution derived from the pETNNs closely approximates the high-fidelity numerical solution. In Fig.~\ref{Exp3.2b}, we plot the error curves in terms of the vorticity of the flow. Our method reaches competitive accuracy compared with other deep learning methods, e.g., \cite{li2022learning, li2020fourier}, under the same Reynolds number.

\subsection{High-dimensional Heat Equation} \label{sec:heat_eq}

\begin{figure}[ht]
\centering
     \includegraphics[width=0.45\linewidth]{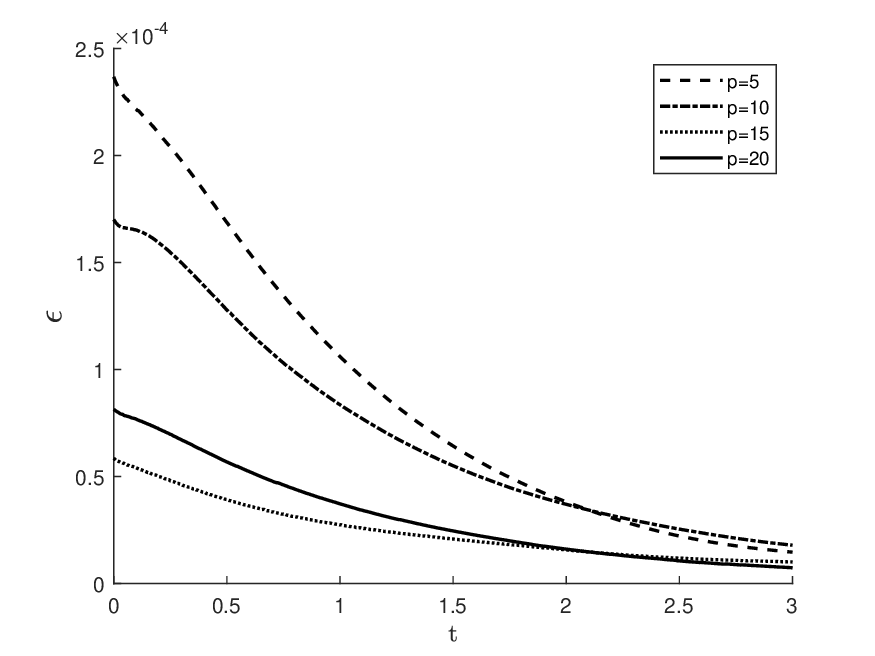}
  \includegraphics[width=0.45\linewidth]{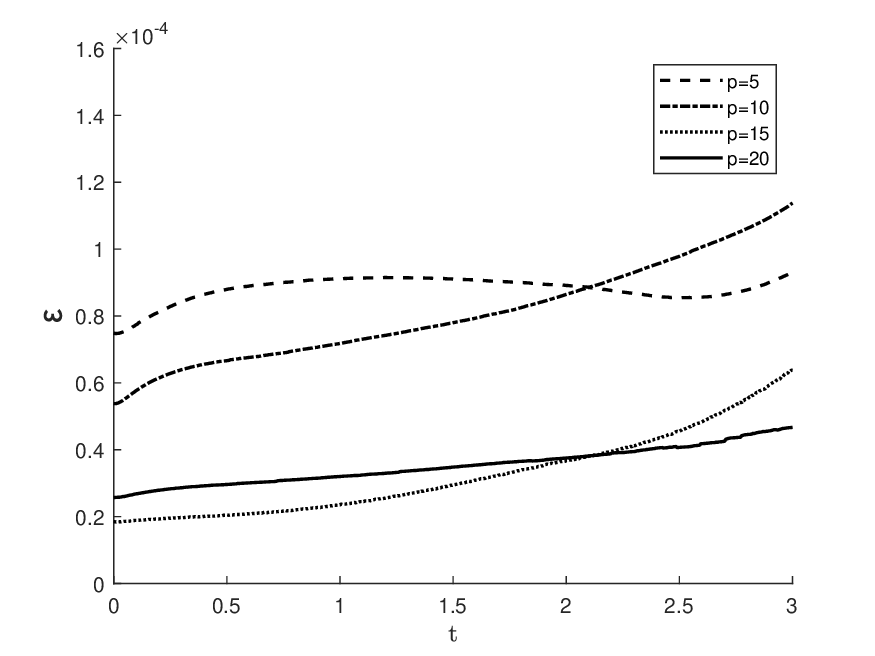}
\caption{The absolute (left) and relative (right) errors of pETNNs for solving 10D heat equation. Here, $p$ denotes the rank of TNNs \eqref{3.1.1} used in pETNNs. We report the error curves in which only $1/3$ of the TNN parameters are updated at each time step.}
\label{Exp3.3}
\end{figure}

The heat equation is a fundamental PDE in mathematics, physics, and engineering that describes the heat distribution in a given region over time. As a parabolic PDE, it is the cornerstone of Fourier theory and heat conduction analysis. Here we consider the 10-dimensional ($d=10$) heat equation,
$$
\frac{\partial u}{\partial t}=v \Delta u, \quad x \in \Omega=[-1,1]^{d}, \; v=\left( (d+3) \pi^2\right)^{-1},
$$
with the initial values
\begin{equation}\label{eq:heat_ini}
u(x, t=0)=\sum_{k=1}^{d} \sin \left( 2\pi x_k \right)\cdot\prod_{i\neq k}^{d} \sin \left(\pi x_i\right),
\end{equation}
and the Dirichlet boundary conditions and its analytical solution is
\begin{equation} \label{eq:heat_sol}
u(x,t)=\sum_{k=1}^{d} \sin \left( 2\pi x_k \right)\cdot\prod_{i\neq k}^{d} \sin \left(\pi x_i\right)\exp (-t).
\end{equation}
In this experiment, each sub-network has three hidden layers with 20 neurons in each layer. We use an RK4 scheme with a fixed time step size of $5e-3$. {   The runtime of pETNNs to reach three times units are 4,225s, 6,088s, 7,580s, and 11,183s for $p = 5,\; 10,\; 15,$ and $20$, respectively. Fig.~\ref{Exp3.3} illustrates the error curves of the numerical solutions obtained from pETNNs, indicating their high accuracy and demonstrating the effectiveness of pETNNs.} In particular, notice that when $p=5$ the rank of TNN is less than the dimension $d=10$, but pETNNs still produce decent approximations. Such a phenomenon suggests the exact solution \eqref{eq:heat_sol} may yield a low-rank structure. We further estimate the rank of the solution tensor and compute its CP decomposition. The corresponding numerical results are reported in Appendix~\ref{app:rank}, which support the existence of a low-rank approximation to the exact solution.

\subsection{High-dimensional Transport Equation} \label{sec:hdtran}

\begin{figure}[ht]
  \centering
  \includegraphics[width=0.45\textwidth]{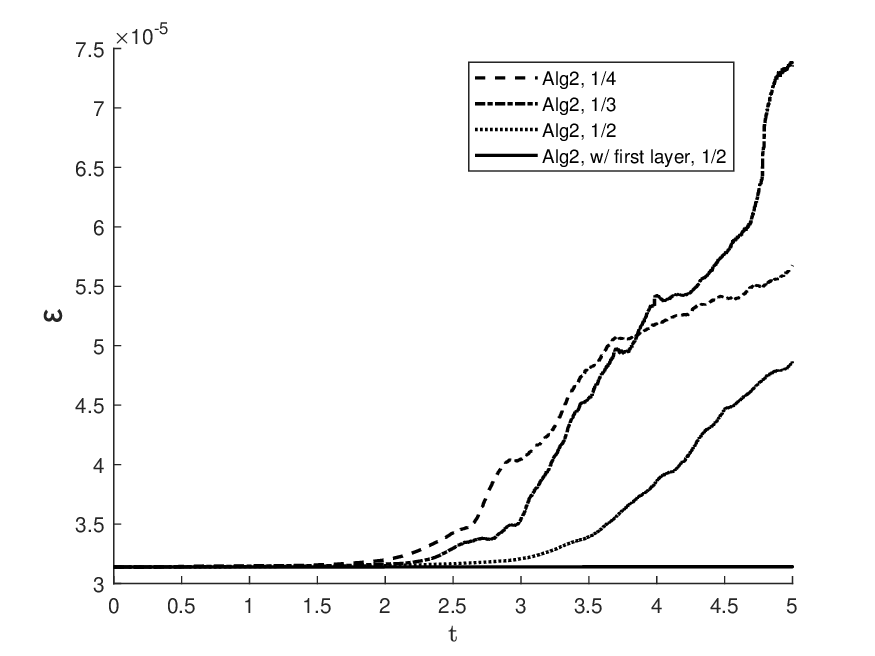}
  \includegraphics[width=0.45\textwidth]{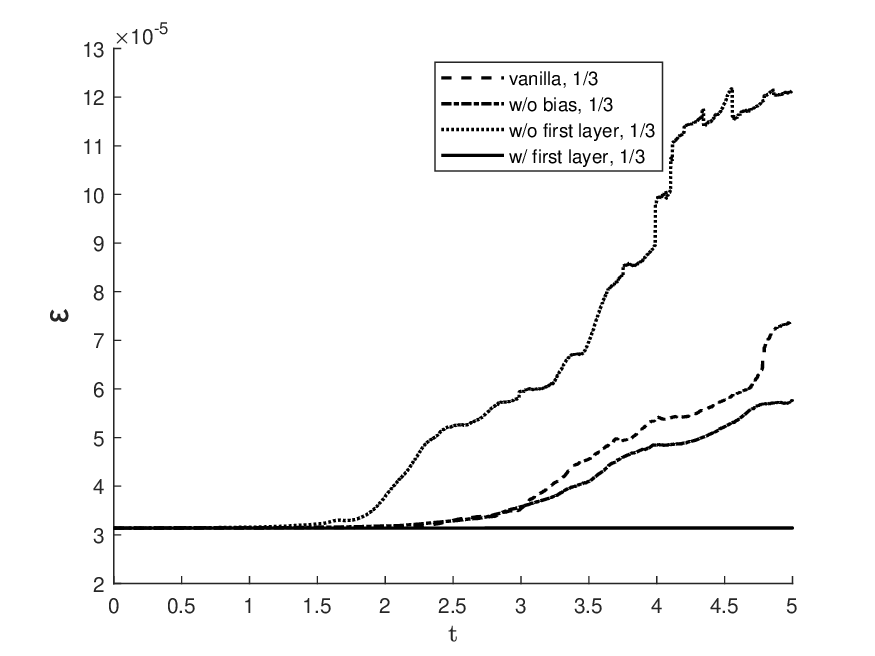}
  \caption{Relative errors of pETNNs with randomized parameter update strategy for solving the 10D transport equation: changing the size of evolutionary parameters (left) and changing the selecting strategy of evolutionary parameters (right). In the left figure, $1/4$, $1/3$, and $1/2$ denote that $1/4$, $1/3 $, and $1/2$ of parameters are updated at each time step. In the right figure, $``vanilla"$ denotes that the update parameters are randomly chosen in every time step;  $``w/\;first\;layer"$ denotes that the inclusion of the update parameters in the first layer is mandatory, whereas $``w/o\;first\;layer"$ signifies precisely the opposite, indicating the exclusion; and $``w/o\;bias"$ denotes the update parameters without the TNNs bias parameter in every time step.}
  \label{Exp3.4a}
\end{figure}

We have used 2D and 3D transport equations with isotropic velocities to validate our pETNNs in Section~\ref{sec:compare}. Here, we will adopt the pETNNs to solve the 10D transport equation to demonstrate the efficacy of our method for solving high-dimensional transport PDE with anisotropic velocities. Besides the increase in dimension, we will also consider more complex wave speeds and initial conditions 
$$
\frac{\partial u}{\partial t}+ \sum_{i=1}^d c_i \frac{\partial u}{\partial x_i}=0, \quad x \in \Omega=[-1,1]^d, \quad u(x, t=0)= \sum_{k=1}^{d} \cos\left( \pi x_k \right) \prod_{i\neq k}^d \sin \left(\pi x_i\right),
$$
with periodic boundary conditions. The analytical solutions are
$$
u(x, t)=\sum_{k=1}^{d} \cos\left( \pi (x_k- c_k t) \right) \prod_{i\neq k}^d \sin \left(\pi (x_i - c_i t)\right),
$$
with $d=10$ and the wave speed $c_{k} = 0.2*k$. We set that each sub-network consists of two hidden layers and each layer has 20 neurons. We use an RK4 scheme with a fixed time step size of $1e-3$. {  The runtime to reach five time units of pETNNs with updating $1/4, 1/3$ and $1/2$ of the parameters are 10,853s, 20,295s, and 59,729s, respectively.} Fig.~\ref{Exp3.4a} shows the results of employing Algorithm~\ref{algorithm2} to the 10D transport equation. The small relative errors demonstrate that our pETNNs perform well in solving high-dimensional transport equations.

Additionally, we have also investigated the impact of parameter selection on numerical accuracy and stability. Within each sub-network, we update $1/3$ of the parameters. Various strategies were employed, including random updates, mandatory updates of the first layer's parameters with each iteration, deliberately excluding the first layer's parameters from the update, and omitting bias from the updating process. The numerical results of these strategies are illustrated in Fig.\ref{Exp3.4a}. The experimental results indicate little difference in the precision of the results obtained by these update methods, which underscores the efficiency of Algorithm \ref{algorithm2} in the short-time regime.

\subsection{A High-dimensional Simplified KdV Equation}\label{sec:KdV}

\begin{figure}[ht]
\centering
  \includegraphics[width=0.45\linewidth]{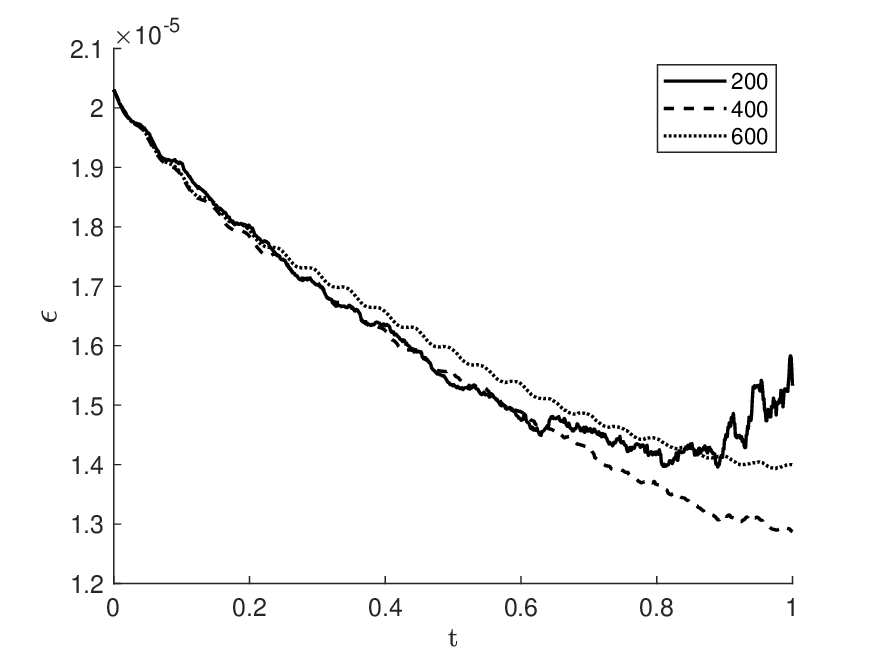}
  \includegraphics[width=0.45\linewidth]{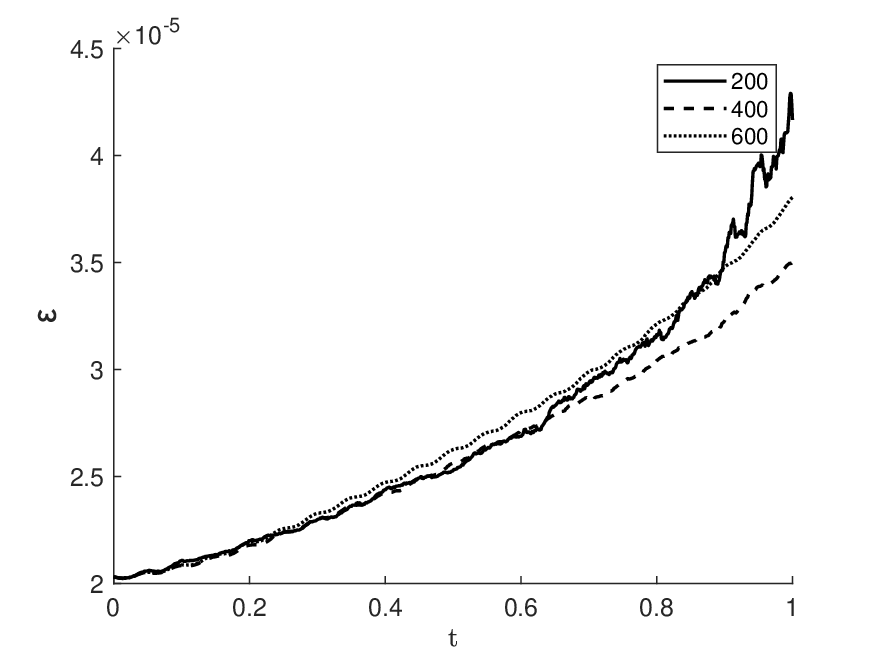}
\caption{Absolute (left) and relative (right) errors of pETNNs for solving the 10D KdV-type equation. Here 200, 400, and 600 parameters are sampled from all of the 1330 parameters in each sub-network and updated in every time step.}
\label{Exp3.5a}
\end{figure}

\begin{figure}[ht]
\centering
  \includegraphics[width=0.45\linewidth]{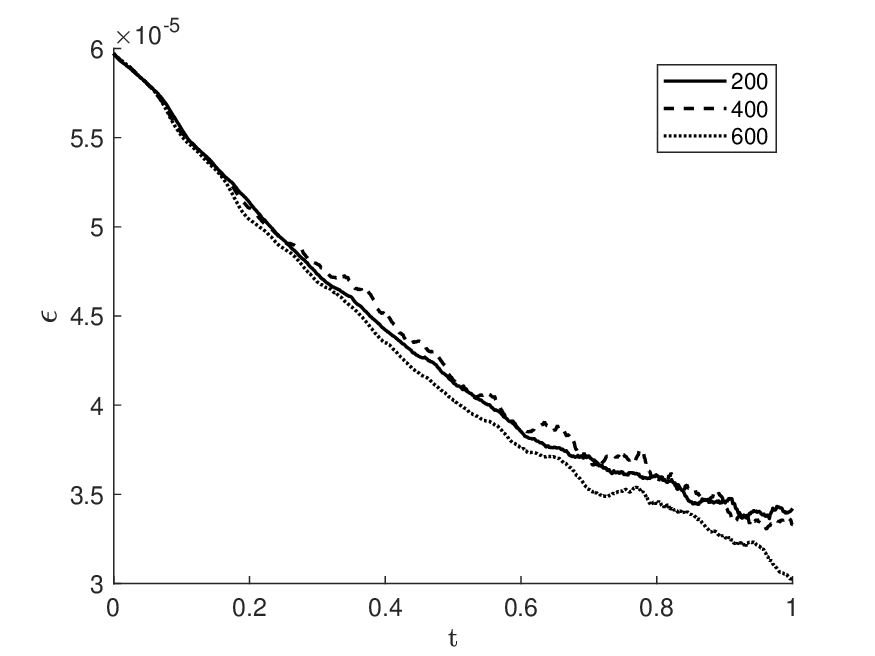}
  \includegraphics[width=0.45\linewidth]{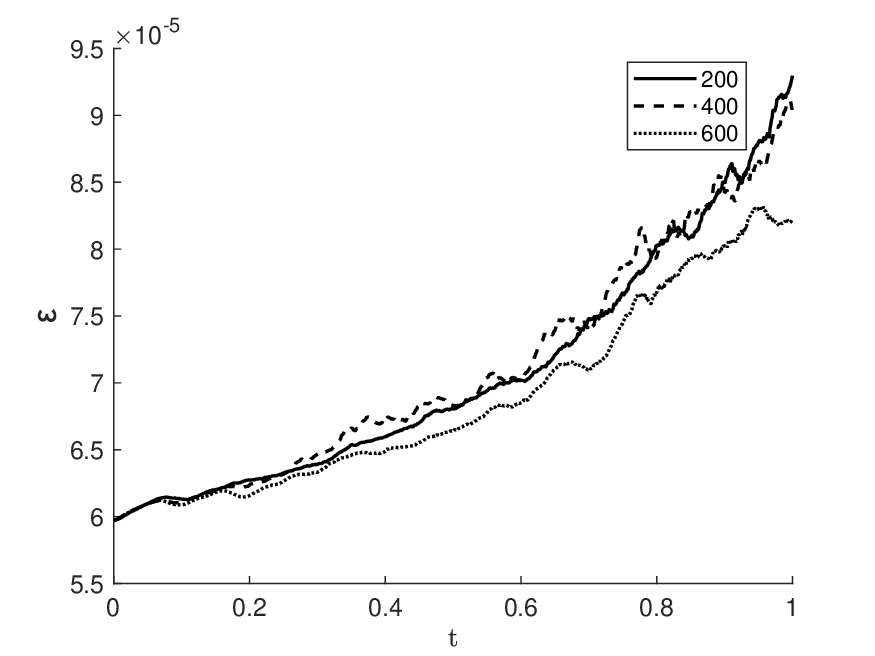}
\caption{Absolute (left) and relative (right) errors of pETNNs for solving the 20D KdV-type equation. Here 200, 400, and 600 parameters are sampled from all of the 1640 parameters in each sub-network and updated in every time step.}
\label{Exp3.5c}
\end{figure}

Our last experiment aims to highlight the ability of pETNNs to handle high-dimensional PDEs involving higher-order derivatives. We will consider the following linear dispersive equation, which can be interpreted as a simplification of the KdV equation,
\begin{equation*}
\frac{\partial u}{\partial t}+c \sum_{i=1}^d \frac{\partial^3 u}{\partial x_i^3}=f, \quad x \in \Omega=[-1,1]^d, \; c=1 / \pi^3.
\end{equation*}
We propose the exact solution 
$$
u(x, t)=\prod_{i=1}^d \sin \left(\pi x_i\right) \exp (-t)
$$
with the corresponding source term
$$
f=-\prod_{i=1}^d \sin (\pi x_i) \exp (-t)-\sum_{i=1}^d \cos (\pi x_i) \prod_{j=1, j \neq i}^d \sin (\pi x_j) \exp (-t).
$$
Here we choose $d=10$ and $20$, and for all cases, the sub-network has two hidden layers with 30 neurons in each layer. The absolute and relative errors are depicted in Figs.~\ref{Exp3.5a} and \ref{Exp3.5c}. All results support pETNNs' excellent performances  for high-dimensional simplified KdV-type equations with high accuracy.

\section{Conclusions and Discussions}

The paper introduces the pETNNs framework as a mesh-free deep learning method for solving time-dependent PDEs. The pETNNs blend the structural advantages of TNNs for high-dimensional problems with evolutionary parameters, yielding not only superior accuracy compared to EDNNs but also exceptional extrapolation capabilities. By employing partial parameter updating strategies, the proposed pETNNs can reduce computational costs significantly. We demonstrate their computational efficiency in addressing high-dimensional time-dependent problems, expanding the utility of deep learning beyond traditional applications in PDEs.

Numerical experiments conducted on various challenging equations, such as the incompressible Navier-Stokes equations, high-dimensional heat equations, high-dimensional transport equations, and dispersive equations of higher-order derivatives, highlight the superior performance of pETNNs in accurately solving these complex problems. {   Although, except for the incompressible Navier-Stokes equations, the PDEs considered yield exact solutions that are separable in space, which can be approximated by TNNs of relatively small size, the separable structure does not alleviate the challenges, e.g., dimensionality, in the numerical implementations. Since pETNNs provides evolutionary parametric approximations by solving the ODE systems derived by the Dirac-Frenkel variational principle, instead of directly fitting the exact solution. Therefore, such test examples allow us to explore the computational efficiency of pETNNs in addressing high-dimensional time-dependent problems based on shallow TNNs, ensuring that computational costs remain manageable.} The findings indicate that pETNNs make a noteworthy contribution to the field of computational methods for time-dependent equations and could play a significant role in advancing scientific research across various disciplines.

{ 
The numerical experiments also reveal several computational challenges. For example, the linear system in \eqref{3.2.4} is often singular, and we have to solve it in the least-squares sense. In general, the singularity of the system is closed related to the PDE formulations, e.g., the variational form and the network structures to parameterize the solution, which will be addressed in future work. In addition, to use the efficient Gaussian quadrature formula, all test problems share rectangular domains, which limits the applicability of pETNNs. As for computational complexity, while the parameter size of a TNN increases linearly with dimension and rank (see the caption of Fig.~\ref{fig1} for details), the runtime for a single iteration increases superlinearly as the time-variant parameter grows.}

To summarize, although the implementation of pETNNs has demonstrated success, future research will aim to comprehensively investigate the scalability and applicability of the method. {  This investigation of scalability will include the adaptation of pETNNs to accommodate irregular domains and more complex boundary conditions,} a systematic examination of parameter update methodologies and an evaluation of long-term stability. In terms of applicability, efforts will be directed towards employing this method to specific physical models, notably in addressing challenges raised in Bose-Einstein condensates \cite{yin2024revealing}, quasicrystals \cite{pnas.2106230118}, and liquid crystals \cite{PhysRevLett.124.090601}, etc.

\section*{Code and Data Availability}
Data and code are available from the corresponding author upon reasonable request for research purposes.

\section*{Acknowledgements}
This work was supported by the National Natural Science Foundation of China (No.12225102, T2321001, 12288101, and 12301520) and the Beijing Outstanding Young Scientist Program \\ (JWZQ20240101027).

\appendix

\section{The numerical integration of the coefficient matrix} 
\label{app:numerical integration}

In Section~\ref{sec:pETNN}, we introduced the governing equation \eqref{3.2.3} of the TNN parameters, where the velocity corresponds to the solution to the following least squares problem
\begin{equation}\label{b.1}
\min J(\gamma) = \frac{1}{2}\int_{\Omega}\left|\frac{\partial {u}}{\partial \theta}\gamma-\mathcal{N}({u})\right|^2dx. \nonumber
\end{equation}
(Here, we have removed the ``hat'' notation for simplicity) To obtain the first-order optimality condition in Eq.~\eqref{3.2.4}, previous works, e.g., \cite{NEURIPS2023_0cb310ed, PhysRevE.104.045303, kast2023positional}, selected a set of spatial points, $\bm{x}_1, \ldots, \bm{x}_n \in \Omega$, and then assemble the observation matrix of the gradients, 
\begin{align}\nonumber
  \textbf{J} = \left[\nabla_{\theta} {u}(\bm{x}_1; \theta), \ldots, \nabla_{\theta} {u}(\bm{x}_n;\theta) \right]^T,
\end{align}
where $\nabla_{\theta} {u} = \left(\frac{\partial {u}}{\partial \theta}\right)^{\top}$ is a column vector. When $n$ is large enough, the ensemble average $\frac1n \textbf{J}^T\textbf{J}$ is a good approximation of the integral $\int_{\Omega}\left(\frac{\partial {u}}{\partial \theta}\right)^{\top}\left(\frac{\partial {u}}{\partial \theta}\right) d x$.

In comparison, thanks to the tensor expansion structure of TNN, the coefficients in the optimality condition can be directly calculated via a numerical integration scheme and this computation is highly parallelizable\cite{wang2023tensor}. For example, suppose we have two rank-p TNNs
\begin{align}\nonumber
u ( x )=\sum_{j=1}^{p} \prod_{i=1}^{d} u_{i, j} ( x_{i} ), \quad v ( x )=\sum_{j=1}^{p} \prod_{i=1}^{d} v_{i, j} (x_{i})
\end{align}
with input $x=(x_1,\ldots,x_d) \in \Omega=[0,1]^{d}$. In each one-dimensional domain, we choose $K$ Gauss points $\{z_k\}_{k=1}^K$ and the corresponding weights $\{w_k\}_{k=1}^K$, then the integration of the product of $u$ and $v$ can be calculated via
\begin{equation*}
\begin{aligned}
\int_{\Omega} u v d x & =\sum_{1 \leq j_1, j_2 \leq p} \int_{\Omega} \prod_{i=1}^d u_{i, j_1} v_{i, j_2} d x =\sum_{1 \leq j_1, j_2 \leq p}\left(\int_0^{1} u_{1, j_1} v_{1, j_2} d x_1\right) \cdots\left(\int_0^{1} u_{d, j_1} v_{d, j_2} d x_d\right) \\
& \approx \sum_{1 \leq j_1, j_2 \leq p}\left(\sum_{k=1}^K w_k u_{1, j_1}\left(z_k\right) v_{1, j_2}\left(z_k\right)\right) \cdots\left(\sum_{k=1}^K w_k u_{d, j_1}\left(z_k\right) v_{d, j_2}\left(z_k\right)\right).
\end{aligned}
\end{equation*}
As for the coefficients in the first-order optimality conditions \eqref{3.2.4}, note that the parameters $\theta$ of TNN can be divided into disjoint components $\theta^1, \ldots, \theta^d$ as in \eqref{3.1.1}, where $\theta^i$ corresponds to the parameters of the sub-network associated with the i-th dimension of the input. As a result, ${u}$ in \eqref{3.2.4} is of the form
\begin{equation*}
    u ( x, \theta )=\sum_{j=1}^{p} \prod_{i=1}^{d} {u}_{i, j} ( x_{i}, \theta^{i} ).
\end{equation*}
For any two parameters $\theta_s$ and $\theta_t$ in $\theta$, we have the following two cases.

\textbf{Case I}: (diagonal blocks) $\theta_s$ and $\theta_t$ are located in the same sub-network. Without loss of generality, we assume that $\theta_s, \theta_t \in \theta^{i_0}$, then the $(s,t)$-component of the matrix $\int_{\Omega}\left(\frac{\partial {u}}{\partial \theta}\right)^T\left(\frac{\partial  {u}}{\partial \theta}\right) d x$ can be calculated by 
\begin{equation*}
\begin{aligned}
\int_{\Omega} \frac{\partial  {u}}{\partial \theta_s} \frac{\partial  {u}}{\partial \theta_t} d x
&\approx
\sum_{1 \leq j_{1}, j_{2} \leq p} \left[ \prod_{i \neq i_0} \left( \sum_{k=1}^{K} w_{k} u_{i , j_{1}} ( z_{k} ) u_{i , j_{2}} ( z_{k} ) \right) \right] \left( \sum_{k=1}^{K} w_{k} \frac{\partial u_{i_0, j_{1}} ( z_{k} ) } {\partial\theta_{s} }\frac{\partial u_{i_0, j_{2}} ( z_{k} )} {\partial\theta_{t}}  \right),
\end{aligned}
\end{equation*}
where the factor $\prod_{i \neq i_0} \left( \sum_{k=1}^{K} w_{k} u_{i , j_{1}} ( z_{k} ) u_{i , j_{2}} ( z_{k} ) \right)$ is shared by all $\theta_s, \theta_t \in \theta^{i_0}$.

\textbf{Case II}: (off-diagonal blocks) $\theta_s$ and $\theta_t$ are located in different sub-networks. We assume that $\theta_s\in \theta^{i_1}$ and $\theta_t \in \theta^{i_2}$, then we have
\begin{equation*}
\begin{aligned}
\int_{\Omega} \frac{\partial  {u}}{\partial \theta_s} \frac{\partial  {u}}{\partial \theta_t} d x
\approx & \sum_{1 \leq j_1, j_2 \leq p}\left[\prod_{i \neq i_1, i_2}\left(\sum_{k=1}^K w_k u_{i , j_1}\left(z_k\right) u_{i , j_2}\left(z_k\right)\right)\right] \\
& \left(\sum_{k=1}^K w_k \frac{\partial u_{i_1, j_1}\left(z_k\right) }{\partial \theta_{s}}u_{i_1, j_2}\left(z_k\right)\right)\left(\sum_{k=1}^K w_k u_{i_2, j_1}\left(z_k\right) \frac{\partial u_{i_2, j_2}\left(z_k\right)}{\partial \theta_{t}}\right),
\end{aligned}
\end{equation*}
where the factor $\prod_{i \neq i_1, i_2}\left(\sum_{k=1}^K w_k u_{i , j_1}\left(z_k\right) u_{i , j_2}\left(z_k\right)\right)$ is shared by all $\theta_s\in \theta^{i_1}$ and $\theta_t \in \theta^{i_2}$.

In practice, we treat the coefficient matrix $\int_{\Omega}\left(\frac{\partial {u}}{\partial \theta}\right)^T\left(\frac{\partial  {u}}{\partial \theta}\right) d x$ as a block matrix of $d\times d$ blocks, and each block can be parallelly computed in a matrix-wise manner, e.g.,  using the ``einsum'' function in {\sf NumPy}.

\section{Tensor rank of the solution to the Heat Equation} \label{app:rank}

In Section~\ref{sec:heat_eq}, Fig.~\ref{Exp3.3} suggests a low-rank TNN ($p=5$) can decently approximate the solution of the heat equation, which reveals the possible existence of a low-rank structure of the solution. In this appendix, we review the basic concept of tensor rank and numerically estimate the rank of the solution to the heat equation in Section~\ref{sec:heat_eq}. For the details of the tensor theory and tensor decompositions, we refer the readers to \cite{2009Tensor, Markus2016Tensor}. We used the {\sf MATLAB} package Tensorlab\cite{TensorLab} to perform the tensor computations within the appendix. A PDF version of the user guide can be found here\footnote{https://www.tensorlab.net/userguide3.pdf}.

An $N$-way tensor $\bm{X} = (x_{i_1i_2\cdots i_N})\in \mathbb{R}^{I_1 \times I_2 \times \cdots \times I_{N}}$ is rank one if it can be written as the outer product of $N$ vectors, i.e.,
\begin{equation*}
    \bm{X} = \bm{a}^{(1)} \circ \bm{a}^{(2)} \circ \cdots \circ \bm{a}^{(N)}, \quad \bm{a}^{k} \in \mathbb{R}^{I_{k}}, \quad k = 1,2,\dots, N.
\end{equation*}
In other words, each element of $\bm{X}$ is the product of the corresponding vector elements:
\begin{equation*}
    x_{i_1i_2\cdots i_N} = a_{i_1}^{(1)} a_{i_{2}}^{(2)} \cdots a_{i_{N}}^{(N)}, \quad 1\leq i_{k} \leq I_{k}, \quad k = 1,2,\dots, N.
\end{equation*}
The CP decomposition factorizes a tensor into a sum of rank-one tensors. For example, given a $3$-way tensor $\bm{Y} \in \mathbb{R}^{I\times J\times K}$, we wish to write $\bm{Y}$ as
\begin{equation*}
    \bm{Y} \approx \sum_{r=1}^{R} \bm{a}_{r}\circ \bm{b}_{r} \circ \bm{c}_{r}, 
\end{equation*}
where $\bm{a}_{r}\in \mathbb{R}^{I}, \bm{b}_{r} \in \mathbb{R}^{J},$ and $\bm{c}_{r} \in \mathbb{R}^{K}$ for $r = 1,2,\dots, R$, and $R$ stands for the rank of the decomposition. In general, the rank of a tensor $\bm{X}$, denoted $\mathrm{rank}(\bm{X})$, is defined as the smallest number of rank-one tensors in an exact CP decomposition. The definition of tensor rank is analogous to the definition of a matrix rank, but their properties are quite different, and there is no straightforward algorithm to determine the rank of a given tensor\cite{2009Tensor}.

Consider the heat equation discussed in Section~\ref{sec:heat_eq}, and introduce its $d$-way solution tensor $\bm{u}(t): [0, T] \rightarrow \mathbb{R}^{n \times n \times \dots \times n}$ based on a homogeneous Cartesian grid over the box $[0,1]^{d}$. The spatial discretization is homogeneous in the sense that the sets of grid points are identical in each direction. As a side note, any straightforward spacial discretization of high-dimensional PDEs inevitably leads to the curse of dimensionality. To address the exponential growth of the computational cost and the storage requirement, a series of techniques have been proposed, e.g., the hierarchical Tucker tensor format. We refer the readers to \cite{Daniele1} and the references therein for the details. Here, since our main interest is the rank of the solution, we consider a Cartesian grid as the spacial discretization for simplicity.

\begin{figure}[ht!]
  \centering
\includegraphics[width=0.45\textwidth]{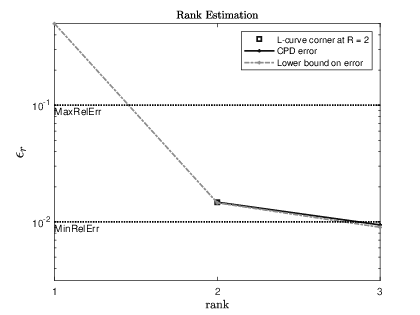}
\includegraphics[width=0.45\textwidth]{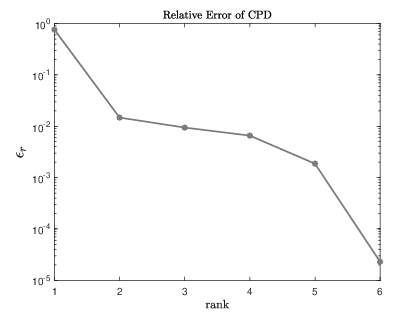}
  \caption{The rank estimation of the solution tensor $\bm{u}(0)$ (left) suggests the lower bound of the rank is 2. The relative error of CP decompositions (right) reproduces the ``L-curve'' as the rank of the decomposition increases. The relative error is defined based on the Frobenius norm.}
  \label{fig:rank}
\end{figure}

Due to the variable separation form of the explicit solution in \eqref{eq:heat_sol}, we know that $\mathrm{rank}(\bm{u}(t))$ is independent of $t$. Set $t=0$, we rewrite the initial condition in \eqref{eq:heat_ini} as
\begin{equation*}
    u(x, t=0) = 2\sum_{k=1}^{d}\cos(2\pi x_{k}) \cdot \prod_{i=1}^{d} \sin (\pi x_{i}),
\end{equation*}
which suggests $\mathrm{rank}(\bm{u}(0)) \leq d$. In particular, we implemented ``rankest'' and ``cpd'' in Tensorlab \cite{TensorLab} to estimate the lower bound of $\mathrm{rank}(\bm{u}(0))$ and the residual error of its CP decomposition, respectively. Fig.~\ref{fig:rank} shows the corresponding numerical results. Due to the storage issue, we set $d=6$ and $n=10$, and we added a $0.5\% $ random noise tensor to the solution tensor $\bm{u}(0)$ to improve the numerical stability. Although we cannot confirm the tensor rank of $\bm{u}(0)$, Fig.~\ref{fig:rank} supports the existence of a low-rank approximation toward the solution, which partially explains the performance of pETNN in Fig.~\ref{Exp3.3} when $p$ is less than the dimension of the problem.

{ 

\section{Solve 2D Transport Equations using space-time TNNs} \label{app:st_TNN}

In Section~\ref{sec:compare}, we compare pETNNs with various baselines based on 2D transport equations. The accuracy of pETNNs suggests that evolutionary TNNs \eqref{3.2.2} form a suitable hypothesis space for 2D transport equations. Thus, a natural question is whether we can mimic the time-marching strategy\cite{wight2020solving} for PINNs to split the time domain and train space-time TNNs using the loss functions.

By treating the time variable as another input, the space-time TNNs propose to approximate the solution to a time-dependent PDE by
\begin{equation*}
 u(x,t) \approx \hat{u}(x,t; \theta) = \sum_{j=1}^p \hat{u}_{d+1, j}(t,\theta^{d+1} )\prod_{i=1}^d \hat{u}_{i, j} \left(x_i, \theta^{i} \right),
\end{equation*}
where $\hat{u}(x,t; \theta)$ is a rank-$p$ TNN with $d+1$ sub-networks, and the $(d+1)$-th sub-network representing the solution's dependence on the time variable $t$. Here, $\theta^{i}$ denotes the parameter of the $i$-th sub-network, $i=1,2,\dots, d+1$.

The space-time TNNs for the 2D transport equation consists of two hidden layers, each with 30 neurons, and is trained in the same way as the PINNs model in Section~\ref{sec:compare}. For simplicity, we consider the time window $[0,1]$, and test the space-time TNNs over different uniform mesh sizes. None of the outcome outperforms the evolutionary TNNs. Table~\ref{table:time_space_TNN} lists the relative error of space-time TNNs under two extreme time steps. Although the parameters of evolutionary TNNs are functions of time, there are identified by numerically solving the governing ODE systems, which has a step size $\Delta t = 5e-3$. Thus, in terms of the mathematical formulation, the solution obtained by evolutionary TNNs share the same structure as space-time TNNs under the same time step. Thus, similar to PINNs\cite{Krishnapriyan2021}, the failure of space-time TNNs is caused by training rather than the representational ability of TNNs.

\begin{table}[htbp]
\centering
\begin{tabular}{@{}ccccccc@{}}
\toprule
& \multicolumn{6}{c}{relative error $\varepsilon_{r}$}                                             \\ \cmidrule(l){2-7}
     & t=0       & t=0.2       & t=0.4       & t=0.6       & t=0.8       & t=1      \\\cmidrule(l){2-7}
S-T TNN($\Delta t=5E-03$) & 4.45E-04	& 2.93E-03 & 	6.78E-03 & 3.00E-02 & 3.97E-02 & 4.51E-02 \\
S-T TNN($\Delta t=1$) & 2.09E-02	& 0.849 &	0.942 & 0.855 & 0.684 & 0.680  \\
ETNNs($\Delta t=5E-03$) & 1.79E-05   &  2.46E-05  &   4.25E-05    & 6.29E-05   & 8.24E-05  & 9.78E-05
 \\ \bottomrule 
\end{tabular}
\caption{Relative errors of space-time TNN under two extreme time steps $\Delta t = 1$ and $\Delta t = 5e-3$ compared with the result of the evolutionary TNNs (ETNNs). }
\label{table:time_space_TNN}
\end{table}

}

\bibliographystyle{plain}
\bibliography{ref}

\end{document}